\documentclass[a4paper,11pt]{article}
\usepackage{amsmath}
\usepackage{amsthm}
\usepackage{amsfonts}
\usepackage{amssymb}
\usepackage{epsfig}
\usepackage[top=3cm, bottom=2cm, left=2.5cm, right=3cm]{geometry}

\newcommand{\NN}{\mathbb N}

\newtheorem*{acknowledgement*}{Acknowledgements} 
\newtheorem{theorem}{Theorem} 
 
\newtheorem{remark}{Remark}

\title{Bimodality, prion aggregates infectivity and prediction of strain phenomenon} \date{}
\author{
Vincent Calvez
\thanks{{\small D\'epartement de Math\'ematiques et Applications, \'Ecole Normale Sup\'erieure, CNRS UMR8553, 45 rue d'Ulm, F-75230 Paris cedex 05. Email: vcalvez@dma.ens.fr}}
\thanks{\small These authors contributed equally to this work} \and 
Natacha Lenuzza
\thanks{{\small CEA-Institute of Emerging Diseases and Innovative Therapies, Route du Panorama, Bat. 60, F-92265 Fontenay-aux-Roses, Emails: natacha.lenuzza@cea.fr, franck.mouthon@cea.fr, jean-philippe.deslys@cea.fr}} \footnotemark[5] \footnotemark[2]  \and 
Dietmar Oelz
\thanks{{\small Universit\"at Wien, Fakult\"at f\"ur Mathematik, Nordbergstr. 15/C/7, 1090 Wien, Austria. Email: dietmar.oelz@univie.ac.at}} \and Jean-Philippe Deslys \footnotemark[2] \and Pascal Laurent 
\thanks{{\small Ecole Centrale Paris, Laboratoire MAS, Grande Voie des Vignes, F-92290 Ch\^ atenay-Malabry, Email: pascal.laurent@ecp.fr}} \and Franck Mouthon \footnotemark[2] 
\thanks{\small These authors contributed equally to supervise this work} \and Beno\^ \i t Perthame
\thanks{{\small Universit\'e Pierre et Marie Curie-Paris 6, UMR 7598 LJLL, BC187, 4, place Jussieu,  F-75252 Paris cedex 5, and Institut Universitaire de France. Email: perthame@ann.jussieu.fr}}  \footnotemark[6]
}

\begin{document}

\maketitle

\begin{abstract} 
We consider a model for  the polymerization (fragmentation) process involved in infectious prion self-replication and study both its dynamics and non-zero steady state. We address several issues. 
Firstly, we give conditions leading  to size repartitions of PrPsc aggregates that exhibit  bimodal distributions, as indicated by recent experimental studies of prion aggregates distribution \cite{Caughey}. This is achieved by a choice of coefficients in the model that are not constant, thus extending a previous study of the nucleated polymerization model \cite{GreerPW}. Surprisingly, conditions for bimodality do not seem to be the most favourable for prion replication. Secondly,  we show stability results for this steady state for general coefficients where reduction to a system of differential equations is not possible. We use a duality method based on recent ideas developed for population models. These results underline the potential influence of the amyloid precursor production rate in promoting amyloidogenic diseases. Finally, we numerically investigate the influence of different parameters of the model on PrPsc accumulation kinetics, in the aim to study specific features of prion strains. This study suggests that PrPsc aggregate size distribution could be a signature of a strain in a given host and a constraint during the adaption mechanism of the species barrier overcoming, that open experimental perspectives for prion strain investigation.

\end{abstract}

\section{Introduction}

Transmissible spongiform encephalopathies (TSE) are fatal, infectious, neurodegenerative diseases. They include bovine spongiform encephalopathies (BSE) in cattle, scrapie in sheep and Creutzfeldt-Jakob disease (CJD) in human \cite{Aguzzi2006}. The infectious agents responsible for disease transmission, known as prions, present some unusual biological properties (as a high resistance to inactivation by heat or radiation). According to the ``protein-only hypothesis", prions may consist in a misfolded protein (called PrPsc), without any nucleic acid. This hypothesis suggests that PrPsc replicates in a self-propagating process, by converting the normal form of PrP (called PrPc) into PrPsc \cite{Prusiner98}. Many evidences are in favour of an autocatalytic replication of PrPsc, as the generation of infectivity from recombinant proteins \cite{Legname04} or the use of in vitro PrPsc conversion systems, such as the protein misfolding cyclic amplification (PMCA) technique \cite{Castilla}. 

However, the precise mechanism of conversion remains unclear. Moreover, prion infectious agent can exist under different strains, characterized by their incubation period and their lesion profile in brains \cite{Fraser}. In the framework of the protein-only hypothesis, it is supposed that strain diversity is supported by various conformation states of PrPsc, that leads to various biological and biochemical properties \cite{Cobb, Morales}. A critical challenge of prion biology consists in elucidating the mechanism of conversion of PrPc into PrPsc, and therefore how a diversity of strains may exist in the same host (expressing the same PrP molecule). 

To investigate the conversion of PrPc in PrPsc, many relevant mathematical modeling of prion replication have been proposed \cite{Eigen,masel,Mobley,Payne}. Major aim of these modelling is to demonstrate that essential features of prion disease can be explained by purely physico-chemical mechanisms, as supposed by the protein-only hypothesis. In addition, mathematical modeling allows to study the effect of every elementary process in a separate manner \cite{Rubenstein}, what is difficult to do experimentally.

The early proposed model is the heterodimer one. It is based on the conformational change of PrPc into PrPsc after the formation of a heterodimeric complex (PrPc + PrPsc $ \rightarrow $ PrPc*PrPsc $ \rightarrow $  PrPsc*PrPsc $ \rightarrow $ 2 PrPsc). This model does not take into account the aggregation of PrPsc, and thereby fails to explain the association between infectivity and aggregated PrP. Some other mechanisms have been proposed, which are interested in PrP aggregation \cite{masel, Mobley, Cohen, Jarrett, Kulkarni}. Based on fibrilar aggregation, the model which seems by now broadly accepted is the one of nucleated polymerization \cite{masel, GreerPW}. In this approach, PrPsc is considered to be a polymeric form of PrPc. Polymers can lengthen by addition of PrPc monomers, and they can replicate by splitting into smaller fragments. It is worth noting that deterministic \cite{masel, GreerPW,Pochel} as well as stochastic \cite{Rubenstein, Pochel} simulations of this model lead to an unimodal size distribution of PrP aggregates, which seems to be quite insensitive to small variations of parameters \cite{Rubenstein}. Greer et al. \cite{GreDrieWanWeb} recently improved the model and include a mean saturation effect by the whole population of polymers onto the lenghtening process (called general incidence), and polymer joining (through a Smoluchowski coagulation equation). In these models, each aggregate has the same behaviour, regardless to its size. However, recent experimental analysis of relation between infectivity and size distribution of PrPsc aggregates (for PrPsc purified from infected brain \cite{Caughey} or for PrPsc produced by PMCA \cite{Weber1, Weber2}) contradicts this uniform behaviour of PrPsc aggregates. In addition, Silveira et al. shows that bimodal distribution of the polymer size is more likely to occur within the real process \cite{Caughey}.

The goal of this present study is to better understand the unexpected experimental size distribution of prion aggregates in brain, and to investigate the potential implication of this size distribution in the strain phenomenon. To do so, we generalize Masel and coauthors' model so as to take into account aggregate size-dependent parameters.

The paper is organized as follows: in section \ref{sec:model} we recall and review the model of Masel {\em et al.} and its continuous version \cite{GreerPW} which is going to be used. We introduce our main improvement, namely a size-dependent lenghtening factor related to nonuniform infectivity rate; and we discuss the eigenvalue problem which is a key tool to analyze the model. In section \ref{sec:bimodal} we investigate heuristically and numerically the formation of a bimodal distribution of polymers. In section \ref{sec:stability} we prove the stability of the zero steady state in the disease free regime (generalizing partially a result in \cite{GreerPW}). Finally, in section \ref{sec:timescales} we study numerically the influence of different parameters on the dynamics of our model, in a prion strain perspective.

\section{The continuum model}
\label{sec:model}

The following set of coupled differential equations  has been introduced by Masel et al. \cite{masel} in order to model the polymerization (aggregation, fragmentation) process involved in infectious prion self-replication.
It describes the dynamics of the quantity of PrPc, $V(t)$, coupled with the evolution of aggregates of PrPsc which have size $i$, $u_i(t)$, 
\begin{equation}
\left\{
\begin{array}{rcl}
\dfrac{dV}{dt} &=& \displaystyle \lambda - \gamma V - \tau V U + 2 \beta\sum_{i=1}^{n_0-1} \sum_{j>i} i u_j ,
\vspace{.2cm}\\ 
\dfrac{du_i}{dt} &=& \displaystyle - \mu u_i - \beta(i-1) u_i - \tau V (u_i - u_{i-i}) + 2 \beta\sum_{j>i} u_j \ , \quad \mbox{for} \ i\geq n_0 \ . \end{array}\right.  \label{eq:Masel}
\end{equation}
The index $n_0$ denotes the minimal size of PrPsc polymers.  
The quantity $U(t) = \sum u_i(t)$ is the total amount of prion aggregates.
The constant parameters $\lambda, \gamma, \tau, \beta, \mu$ are respectively the basal synthesis rate of PrPc, the degradation rate of PrPc, the conversion rate of PrPc into PrPsc (autocatalytic process following the mass action law), the fragmentation coefficient and the degradation rate of PrPsc.

Analysis is simpler in the framework of continuous size of prions because analytical tools can serve to find simpler formulations. Accordingly, Greer {\em et al.} introduce a continuous version of \eqref{eq:Masel} where the variable $x\in (0, +\infty)$ denotes the size of aggregates and replaces the index $i\in \NN$. Asymptotic derivations of continuous models from discrete models can be found in \cite{ELM, LM}. The continuous model reads, with possibly nonconstant coefficients, 
\begin{equation}
\label{eq:Greer}
\left\{
\begin{array}{rl}
\dfrac{dV(t)}{dt}&=\displaystyle \lambda- V(t)\left[\gamma+ \int_{x_0}^{\infty} \tau(x) u(x,t) \; dx\right] +2 \int_0^{x_0}x \int_{x}^{\infty} \beta(y) \kappa(x,y)\, u(y,t) \, dy \, dx, 
\vspace{.2cm}\\
\dfrac{\partial}{\partial t} u(x,t) &= \displaystyle - V(t) \frac{\partial}{\partial x} \big(\tau(x) u(x,t)\big) - [\mu(x) + \beta(x)] u(x,t) + 2 \int_{x}^{\infty} \beta(y) \kappa (x,y) \, u(y,t) \, dy, 
\vspace{.2cm}\\
u(x_0,t) &=0,
\end{array} \right.
\end{equation}
together with appropriate initial conditions. 
This is a well established family of models used for describing aggregation, fragmentation and possibly coagulation in polymers and also size structured cell dynamics with finite resources  \cite{escobedo, Perthame_LN,michel, michel1}. Well-posedness, in the class of weak solutions,  can be found in \cite{LW, SW}.

The transport term $ V(t) \dfrac{\partial}{\partial x} \big(\tau(x) u(x,t)\big)$ accounts for the growth in size of polymers: their size grows with the speed $V(t) \tau(x) $, proportional to the available PrPc molecules $V(t)$, with an aggregation ability depending on the size of the polymer (a conceivable hypothesis being that their size confers them a peculiar geometry affecting the autocatalytic process). The fragmentation rate, for a polymer of size $y$, is $\beta(y) >0$. The repartition of the two fragments of (smaller) sizes $x$ and $y-x$ is given by $\kappa (x,y) \geq 0$. It should thus satisfy the two usual laws \cite{MD} expressing that the number of fragments increases but with constant total molecular mass (recall the factor $2$ in the right hand side of (\ref{eq:Greer}))
\begin{equation}
\label{as:kappa}
\int_0^y \kappa(x,y) dx = 1\, , \qquad \int_0^y x \, \kappa(x,y) dx = \frac y 2\,  . 
\end{equation}

As usual for size-structured population models \cite{Perthame_LN, MD}, this model may incorporate a minimal size of infectious PrPres aggregates $x_0 \ge 0$, whose value remains unknown. Experimentally, no monomer of PrPsc has been isolated yet. In addition, small aggregates have  been shown not to be infectious \cite{Caughey} (even-though they have to enter the actual modeling) and thus the assumption of a critical size of nucleation $x_0$ has been used. However, the continuous model holds under the assumption that the monomer size can be neglected as opposed to possibly numerous large polymers (such that $x_0\simeq0$). Therefore we keep $x_0$ in the following but we sometimes neglect it in order to simplify the presentation.

\begin{figure}
\centering
\includegraphics[width=.45 \linewidth]{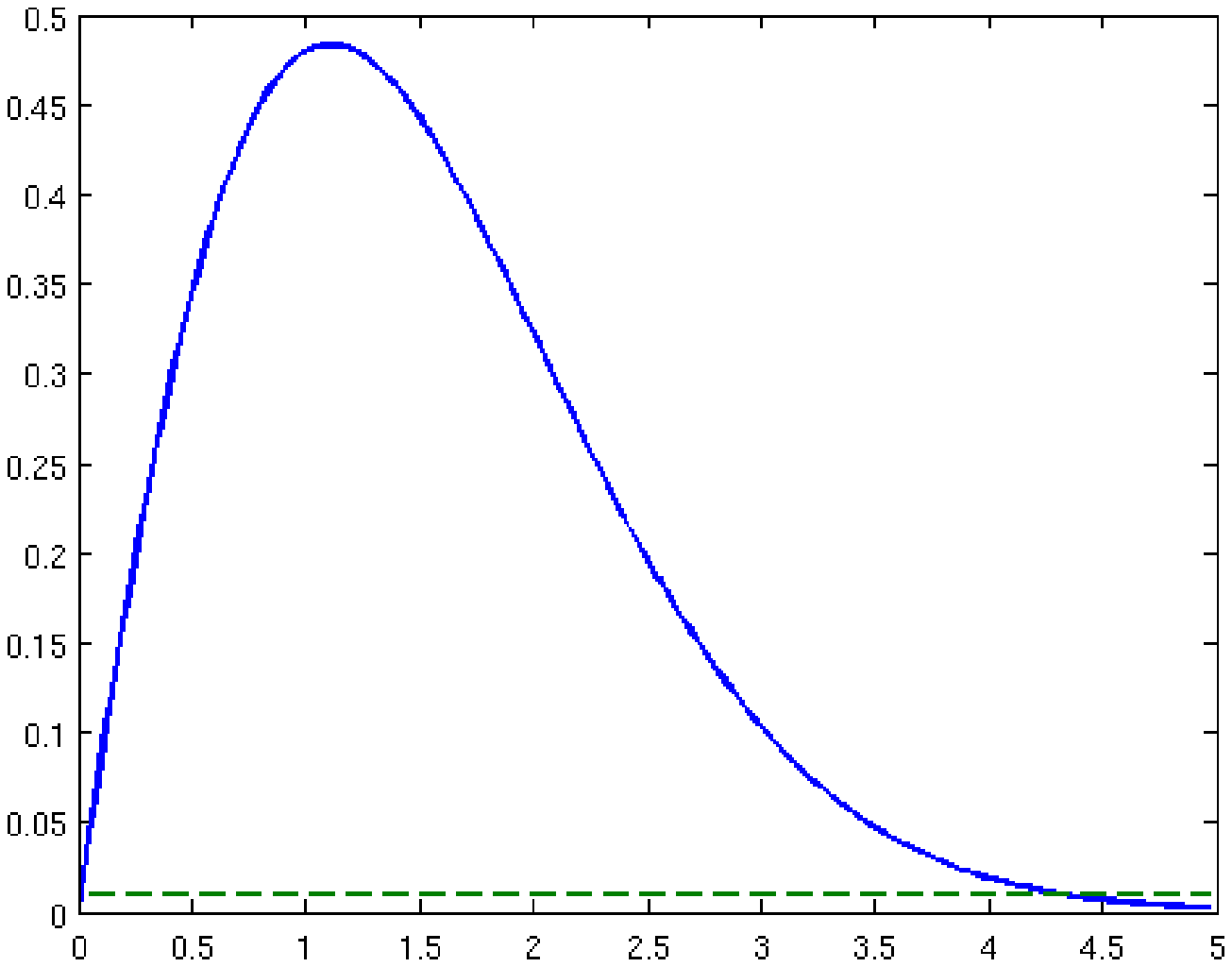}
\includegraphics[width=.45 \linewidth]{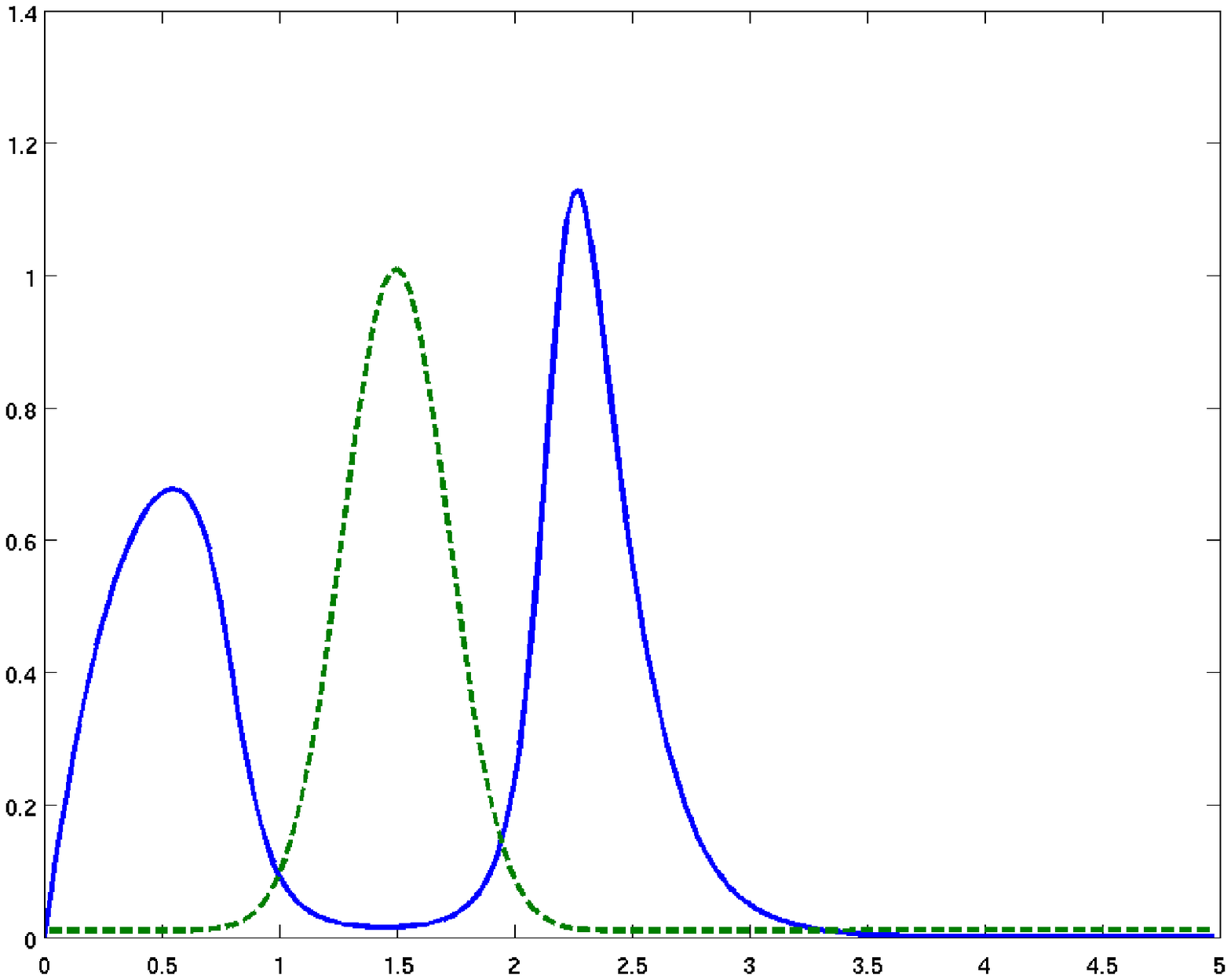}
\caption{
{\sc The splitting effect.} 
(left) A typical stationary distribution for model \eqref{eq:Greer} with 'constant coefficients' according to (\ref{as:cc1}--\ref{as:cc3}). This unimodal shape is given explicitly by a suitable dilation of the profile $\Phi(r) = (r+r^2/2)\exp(-r-r^2/2)$. 
(right)
 For a suitable bell-shaped converting rate $\tau(x)$ (dotted line), the non-zero steady state solution to (\ref{eq:stst}) presents a bimodal shape (solid line). As compared to the constant converting rate case $\tau(x) = \tau_0$ (left), the peak of $\tau$ acts so as to split the distribution into a bimodal one.
\emph{Abscissa = PrPsc aggregates size; Ordinate = $\tau(x)$ either $u_\infty(x)$ in arbitrary units} 
}\label{fig:bimodal} 
\end{figure}

\paragraph{Assumptions on the coefficients.}
In its entire generality, model \eqref{eq:Greer} is rather difficult to attack even though some qualitative behaviors can be described as in section \ref{sec:stability}. We aim to reduce the complexity of this system in order to extract some relevant information as the bimodal distribution motivated by the recent work of Siveira et al. \cite{Caughey}. 
A simple choice for the coefficients 
is as follows \cite{GreerPW}: the fragmentation rate is taken proportional to the fragment size and the degradation rate does not depend of the size, that is 
\begin{equation}
\beta(y) = \beta_0 \; y, \qquad \mu(x)\equiv \mu_0,
\label{as:cc1}
\end{equation}
and the probability distribution of fragments of size $x$ is chosen to be uniform with respect to the length of the splitted polymers of size $y$: 
\begin{equation}
\label{def:kappa}
\kappa(x,y) = \begin{cases}
0 & \text{ if } y \leq x_0 \ \mbox{or}\ y\leq x \\ 		1/y & \text{ if } y > x_0 \ \mbox{and}\ 0<x<y . \end{cases}
\end{equation}
Last, but not least, the transition rate $\tau(x)$ may also be assumed to be constant: 
\begin{equation}
\tau(x)\equiv \tau_0.
\label{as:cc3}
\end{equation} 
In this situation the possible time asymptotic behaviors have been completely classified \cite{masel,GreerPW} as well as the stable distribution of aggregates (see Figure \ref{fig:bimodal}(left) and the paragraph "Constant coefficients`` at the end of this section). When referring to the model of Greer {\em et al.} with a constant rate $\tau_0$, we will denote it by the 'constant coefficients' case. 

With this respect, the main purpose of our work is  to take into account different shapes for the function $\tau(x)$ in order to fit the recent experimental data. Also we perform a stability analysis which does not assume particular coefficients as in (\ref{as:cc1}--\ref{as:cc3}). 
At this stage of knowledge of the biochemical process, the available microscopic experimental data are insufficient to investigate different fragmentation laws in details. Indeed, we are not able to dissociate the three elementary processes (degradation/sequestration, splitting and polymerization) implicated in prion replication. We have made the choice of varying the rate of infection $\tau(x)$ as a first step, whereas Silveira et al. clearly indicate that large polymers are more stable than small ones (this could be assumed as a reduced degradation rate or a reduced fragmentation factor for large sizes).

The system (\ref{eq:Greer}) keeps an important biochemical property, that is the prion molecules are properly transfered from one configuration to another (inducing no loss of mass during fragmentation or polymerization). This enhances the following macroscopic law involving the total quantity of polymers $U$ and the mean length of aggregates $P$: \[
U(t) = \int_{x_0}^{\infty} u(x,t)\ dx \ , \qquad P(t) = \int_{x_0}^\infty x u(x,t)\ dx \ . \]
In fact equation  (\ref{eq:Greer}) yields,
\begin{equation}
\label{eq:conservation}
\frac{d}{dt} \left(V(t) + P(t) \right) =\lambda-\gamma V(t)- \int_{x_0}^\infty \mu( x) x \; u(x,t)\, dx \ . \end{equation}


\begin{figure}
\centering
\includegraphics[width=.6\linewidth]{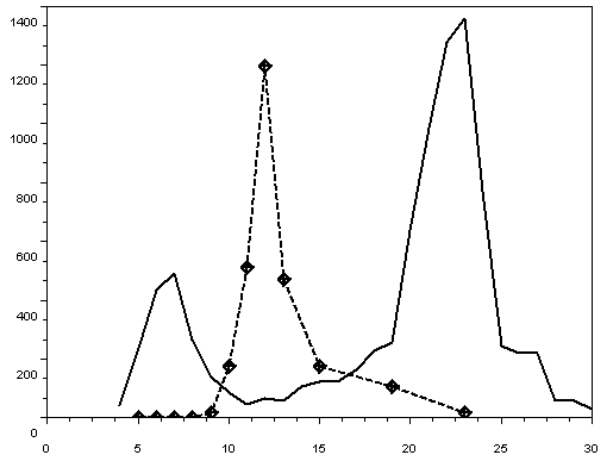} \caption{Experimental analysis of PrPsc aggregates.(Data have been kindly provided by Silveira. For details, see \cite{Caughey}). Size-distribution of PrPsc aggregates  in a whole infected hamster brain (Solid line)  and relative specific infectivity  of each fraction containing PrPsc aggregates with respect to their size (dotted line). The former should correspond to the quantity $xu(x,t)$ within our formalism, and the latter represents a bell-shaped converting rate being $\tau(x)$ in model \eqref{eq:Greer}.
\emph{Abscissa = fraction number; ordinate = PrP (in ng) for the solid line, and relative specific infectivity (in arbitrary units) for the dotted line}. 
Note that the scale of the polymer size $x$ in this graph is not likely to be linear. However we can hardly evaluate the precise scale because of the limitation of measurement experiments. 
}
\label{fig:Silveira}
\end{figure}

\paragraph{Stationary polymer distribution.}
Contrary to the one-peak distribution of infectious prion proteins depicted in Figure \ref{fig:bimodal}(left) under the 'constant coefficients' assumptions (\ref{as:cc1}--\ref{as:cc3}), experiments indicate that this distribution is rather bimodal (Figure \ref{fig:Silveira}).  
It comes out that bimodal configurations could be the result of a non-constant converting factor $\tau(x)$ (having a bell shape for instance).
As we are interested mainly in long-time dynamics resulting in a balance between polymerization and fragmentation, we will often stress out the so-called 'Stable Size Distributions' (see \cite{MD}) that are the stable steady states of system (\ref{eq:Greer}), i.e. satisfying 
\begin{equation}
\label{eq:stst}
\left\{\begin{array}{l}
\displaystyle V_\infty \left[ \gamma + \int_{x_0}^{\infty} \tau(x) u_\infty(x) \; dx\right] =\lambda + 2 \int_0^{x_0}x \int_{x}^{\infty} \beta (y)  \kappa(x,y) \, u_\infty(y) \, dy \, dx , \vspace{.2cm}\\
\displaystyle V_\infty \frac{\partial}{\partial x} \big(\tau(x) u_\infty(x)\big) +[ \mu(x)  + \beta (x)] u_\infty(x) = 2 \int_{x}^{\infty} \beta (y) \kappa(x,y)\, u_\infty(y) \, dy \ ,
\end{array} \right.
\end{equation}
with sufficient decay conditions at infinity. 
There are in general two steady states, the first one corresponds to 'healthy cells', and we refer to it as the 'zero steady state':
\begin{equation}
\label{eq:zerostst}
\overline V= \frac{\lambda}{ \gamma}, \qquad \overline u=0. 
\end{equation}
A second possible steady state, denoted by $(V_\infty,u_\infty)$, corresponds to the infection regime, and can be understood as follows. For a given $V$, the fragmentation equation has a  dominant eigenvalue $\Lambda(V)$ (growth rate), associated to a nonnegative eigenvector. In other words there is a unique solution to 
\begin{equation}
\label{eq:frag}
\left\{\begin{array}{l}
 \displaystyle V \frac{\partial}{\partial x} \big(\tau(x)  {\mathcal U} (V;x)\big) + [ \mu(x) + \beta(x)] {\mathcal U} (V;x) - 2 \int_{x}^{\infty} \beta(y) \kappa (x,y) \, {\mathcal U} (V;y) \, dy=\Lambda(V) \; {\mathcal U} (V;x), 
\vspace{.2cm}\\
{\mathcal U} (V;x_0) =0, \quad {\mathcal U} (V;x) \ge 0, \quad \int {\mathcal U} (V;x)\, dx =1 \ .
\end{array} \right.
\end{equation}
Such solutions have been shown to exist, and the integrability condition implies a much faster decay at infinity, 
see  \cite{Perthame_LN,michel, PeRy}. 

\begin{figure}\centering
\includegraphics[width = .45\linewidth]{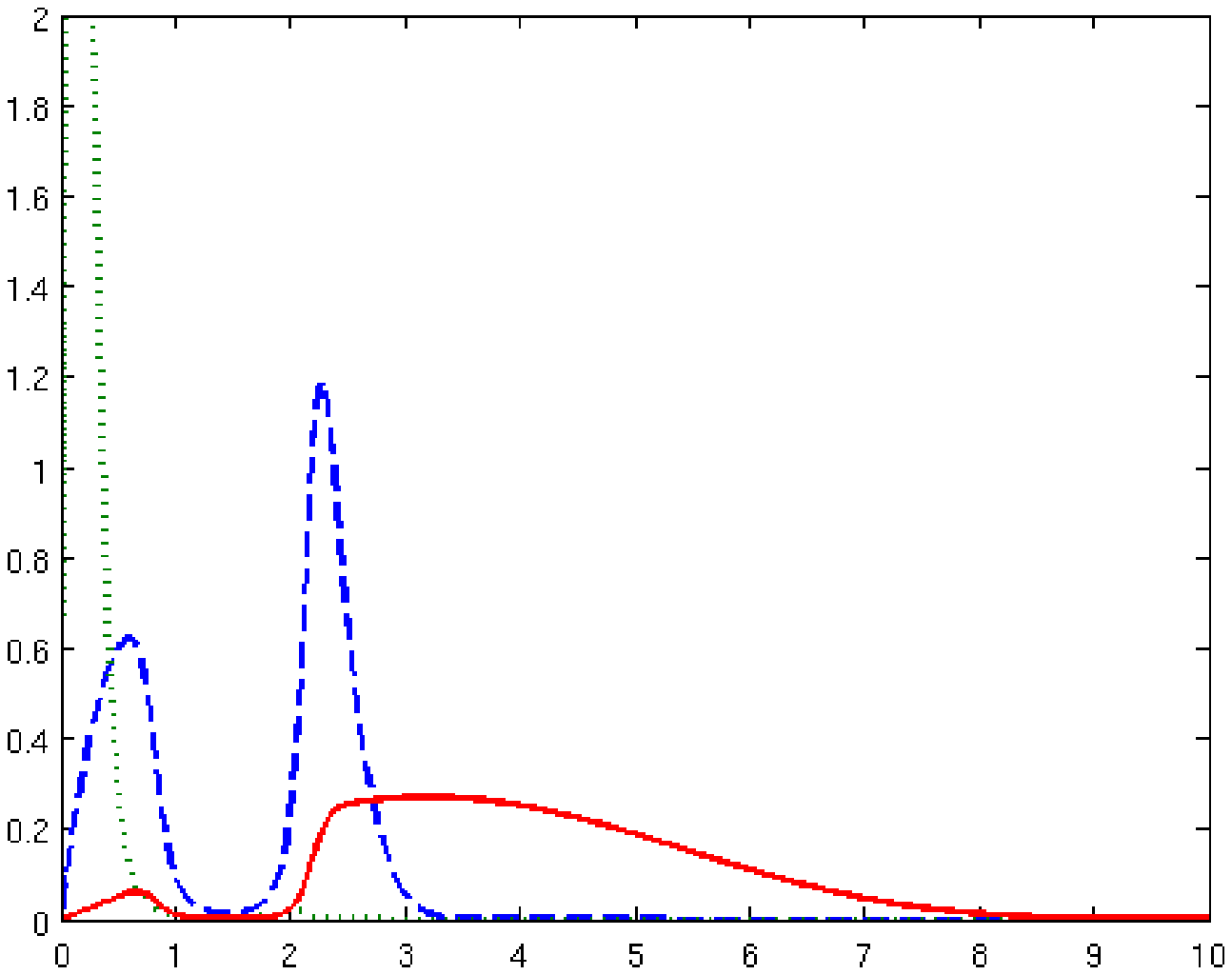}
\includegraphics[width = .45\linewidth]{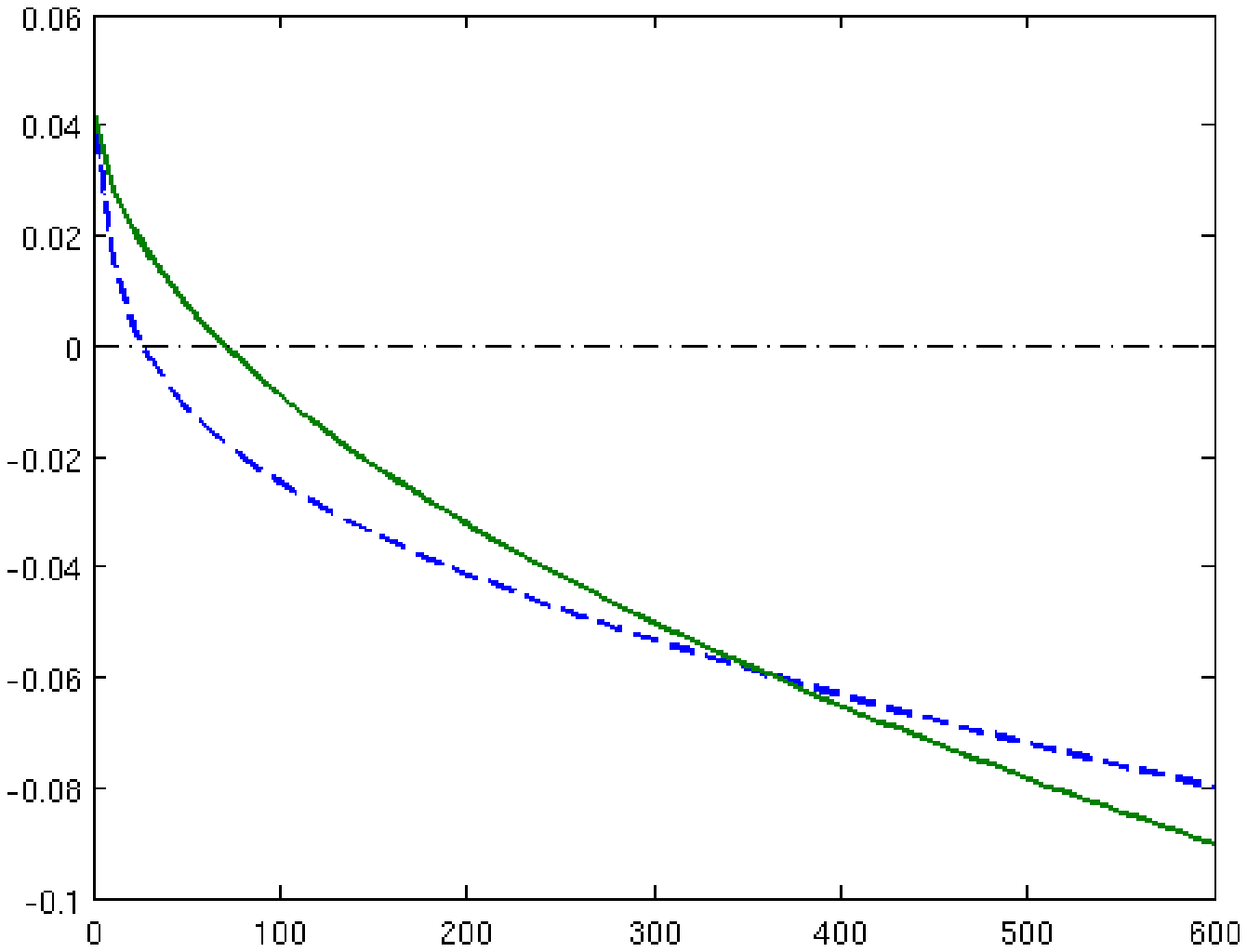}
\caption{{\sc Eigenvalue problem for the microscopic distribution.}
(left) Several eigenfunctions are plotted, for $V$ respectively above (full line), below (dotted line), and close to (dashed line) the equilibrium value $V_\infty$. This is the 'constant coefficient' configuration except for the bell-like converting rate $\tau(x)$. Coefficients' values are given in Section \ref{sec:timescales}.
(right) Numerical computation of the eigenvalue function $\Lambda(V)$ for a bell-like function $\tau(x)$ (dashed line) compared to the 'constant coefficients case', where $\Lambda(V) =  \mu_0- \sqrt{\tau_0 \beta_0 V}$ (full line). 
}
\label{fig:eigenvalue} 
\end{figure}

\begin{remark}
Intermediate asymptotics of \eqref{eq:Greer} and particularly early stages of the infection (exponential phase) are also considered in the sequel (Section \ref{sec:timescales}). Interestingly, the distribution shape corresponding to the exponential phase (Figure \ref{fig:eigenvalue}(left, full line)) matches quantitatively with the experimental data of Figure \ref{fig:Silveira}. Indeed we observe that large polymers are prevalent, corresponding to a level of PrPc close to the healthy state $V\approx \overline V$ (see also Section \ref{sec:timescales}). 
\end{remark}

The equilibrium distribution $u_\infty$ arises as an eigenvector associated to the eigenvalue $\Lambda(V_\infty) = 0$. In fact, this characterizes both the ground level of PrPc, $V_\infty$, and the shape of the polymer distribution $\mathcal U(V_\infty;x)$, up to a constant factor (the total number of polymers). This missing factor is determined thanks to the first equation of \eqref{eq:stst} as we show it later. Interestingly enough, the value $V_\infty$ does not depend on the differential equation driving $V(t)$. In particular it does not depend on $\lambda$ and $\gamma$.

The decay condition at infinity ensures that the following integrations by parts can be justified (we assume $x_0 = 0$ below for the sake of simplicity). 
Integrating \eqref{eq:frag} successively against $1$ and $x$, and using (\ref{as:kappa}), gives respectively
\begin{align*} 
0 +  \int \mu(x) \mathcal U(x)\, dx + \int \beta( x) \mathcal U(x)\, dx - 2  \int_0^{\infty}\int_x^\infty \beta(y) \kappa(x,y)   \mathcal U(y)\, dy &= \Lambda \int \mathcal U(x)\, dx \ , 
\end{align*}
and 
\begin{align*} 
- V\int \tau(x)\mathcal U(x)\, dx + \int (x\mu(x) + x \beta (x)) \mathcal U(x)\, dx -  2  \int_0^{\infty} x \int_x^\infty \kappa(x,y) \beta( y)  \mathcal U (y)\, dy &= \Lambda \int x\mathcal U(x)\, dx \ .
\end{align*}
As a direct consequence we obtain 
\begin{equation}
\label{eq:Lambda}
 \Lambda (V)= \int \big[ \mu(x) -  \beta(x) \big] \mathcal U(V;x)\, dx \frac{- V\int \tau(x)\mathcal U(V;x)\, dx + \int x  \mu(x) \mathcal U(V;x)\, dx}{ \int x\mathcal U(V;x)\, dx}.
\end{equation}

\paragraph{Macroscopic (non zero) steady state.}
Back to \eqref{eq:stst},  the steady state of interest, if it exists, is given according to the rule 
\begin{equation}
\label{eq:stabless}
\Lambda(V_\infty)=0, \qquad u_\infty(x)= \varrho_\infty \; {\mathcal U} (V_\infty;x) \ , 
\end{equation} 
where $\varrho$ denotes the total number of polymers. There do not always exist non-trivial steady states. To see it, consider again the case $x_0=0$. Then, according to the first equation of system \eqref{eq:stst}, we have
\[
V_\infty \left[ \gamma  +  \varrho_\infty \int\tau(x) {\mathcal U} (V_\infty;x)\, dx \right]= \lambda \ ,\]
or, equivalently
\[ \varrho_\infty =\dfrac{ \lambda V_\infty^{-1} - \gamma }{ \int \tau(x) {\mathcal U} (V_\infty;x)\, dx} >0 \ . \]
This points out a constraint for  the non trivial steady state $(V_\infty, u_\infty)$ to exist, namely  it is required that 
\begin{equation}
\label{as:stst}
\gamma V_\infty < \lambda \, , \quad (\mbox{or}\quad V_\infty < \overline V )\,  .
\end{equation}

We tackle in this issue the nonlinear (in)stability of the trivial steady state $(\overline V,0)$ in  section \ref{sec:stability} for general coefficients. In particular, we show that stability of the trivial steady state when $\overline V < V_\infty  $ follows directly from the property that \begin{equation} \Lambda(V)  \quad \mbox{is a decreasing function.}\end{equation}

\paragraph{Constant coefficients.}
The case of 'constant coefficients' (\ref{as:cc1}--\ref{as:cc3}) has been completely understood by 
Greer {\em et al.} \cite{GreerPW}. We have the opportunity to recall their results under the viewpoint of the eigenvalue $\Lambda$ in (\ref{eq:frag}). This eigenvalue can be explicitly computed from \eqref{eq:Lambda} which yields
$$
 \Lambda (V)-\mu_0 = - \beta_0 \int  x \; \mathcal U(V;x)\, dx  \frac{- \tau_0 V }{ \int x \; \mathcal U(V;x)\, dx}.
$$
Eliminating the quantity $\int x\; \mathcal U(V;x)\, dx =  \frac{\tau_0 V}{\beta_0}$, we obtain 
\begin{equation}
\Lambda(V) =  \mu_0  -\sqrt{\tau_0 \beta_0 V} \ .
\label{eq:eigenvaccc}
\end{equation}
 Notice that it is a  decreasing function of $V$ (this property is in fact crucial in the subsequent analysis). We also point out for further use the following consequence for the non-zero steady state:
 \begin{equation}
\frac{ \int x \; u_\infty(x)\, dx}{ \int  u_\infty(x)\, dx}= \frac{\mu_0}{\beta_0}.
\label{eq:peakloc}
 \end{equation}
Following the constraint \eqref{as:stst}, a non trivial steady state exists if and only if $\gamma <\lambda \tau_0$. Still from \cite{GreerPW}, we know that it is globally asymtotically stable. On the contrary, when it does not exist, the zero steady state is  globally asymptotically stable.

\section{Bimodal fragment distribution of PrP} 
\label{sec:bimodal}

Recent experimental data obtained in \cite{Caughey} suggest that PrPsc aggregates have different infectious abilities, according to their size (see Figure \ref{fig:Silveira}). Although infectious ability actually seems to be a competition between the three processes involved in the model (namely fragmentation, polymerization and degradation) one simple way to improve the model compared to \cite{GreerPW} is to introduce a single variable  coefficient. Results from \cite{Caughey} show different in vitro abilities to convert PrPc into PrPsc, arguing in favour of a non-constant polymerization rate. 

Therefore, in this section we consider the system (\ref{eq:stst}) with the choice $\mu(x) = \mu_0$, $\beta(x)=\beta_0 x$ and $\kappa$ given by (\ref{def:kappa}), and a non-constant polymerization rate $\tau=\tau(x)$. We study, by means of numerical experiments and analytical methods, how bimodal fragment distributions may occur. 

\paragraph{Direct numerical simulations.} 

Our first observation is that a reasonable bell shaped converting factor $\tau(x)$ induces a bimodal stationary distribution (see Figure \ref{fig:eigenvalue} for instance). This is supported by various simulations of the model (\ref{eq:stst}). This numerical result leads to investigate analytically which conditions are needed for a bimodal distribution to appear. This is our next step.

\paragraph{A differential equation formulation.}

As a first step, we present a simpler differential formulation of the stationary problem  (\ref{eq:stst}), where we strongly use our particular choice of coefficients (all are constants except $\tau(x)$).

Differentiating equation (\ref{eq:stst}), we find the following second order differential equation
\begin{equation}
\label{eq:ode}
V_\infty\big( \tau( x) u_\infty(x)\big)_{xx}+ \big(  u_\infty(x) \left( \mu_0 + \beta_0 x \right) \big)_x+ 2 \beta_0  u_\infty(x) = 0\ . 
\end{equation}
This differential equation is complemented by the boundary conditions
\[
u_\infty (x_0)=0, \qquad V_\infty (\tau u_\infty )_x\big|_{x=x_{0}}=  2 \beta_0 \int_{x_0}^\infty u_\infty(y) \  .
\]
Recall also the additional relation derived from (\ref{eq:stst}), after testing against $x$,
\begin{equation*}
 V_\infty \int_{x_0}^\infty  \tau(x) u_\infty(x)\, dx =  \mu \int_{x_0}^\infty x \; u_\infty(x)\, dx .
\end{equation*}

The equation \eqref{eq:ode} is reminiscent of the standard one-dimensional Fokker-Planck equation. In fact, if we drop the last term of order 0, it becomes
\begin{displaymath}
\displaystyle{V_\infty} \big(\tau(x) (u_\infty(x))_x\big)_{x} + \left( \psi_x(x) u_\infty(x) \right)_x 
= 0\ ,
\end{displaymath}
that is, a inhomogeneous diffusion with an effective potential given by $\psi(x)={V_\infty} \tau(x) + \mu_0 x  +\frac12 \beta_0 x^2$. Hence the contribution of a non constant $\tau(x)$ might rise  a double well potential on $\psi(x)$, and thus be responsible for the formation of a two-peaked solution. This holds, by perturbation analysis, at least for $\beta_0$ small.  But we can give a more accurate estimate.

\paragraph{Condition for a bimodal equilibrium distribution.} 

We now give a necessary condition for a bimodal distribution. 
We evaluate the derivatives in the microscopic equilibrium equation \eqref{eq:ode}, 
\begin{equation*}
{V_\infty } \big( \tau( x) (u_{\infty})_{xx}( x) +2
\tau_x( x) (u_{\infty})_x( x)+\tau_{xx}( x) u_\infty( x) \big)= \left(- u_\infty \left( \mu_0 + \beta_0 \, x \right) \right)_x- 2 \beta_0 u_\infty(x) . 
\end{equation*}
A bimodal length distribution $u_\infty$ has a critical point where it is convex and therefore we are looking for a point $x^*$ such that $(u_{\infty})_x(x^*)=0$ and $(u_{\infty})_{xx}(x^*) > 0$. A necessary condition on $\tau(x)$ for the existence of such a point is therefore \begin{displaymath} 
\inf_{x > 0}V_\infty \;  \tau_{xx}(x) < -3 \beta_0 .
\end{displaymath}
The above condition is indeed intricated because $V_\infty$ itself depends on $\tau(x)$. But it gives a first  insight of the desirable conditions, meaning that a peak which is concave enough is required.
As we will see it below, this condition is not sufficient. The locus of this strongly concave 'peak' is also important

\begin{figure}
\begin{center}
\includegraphics[width = .45\linewidth]{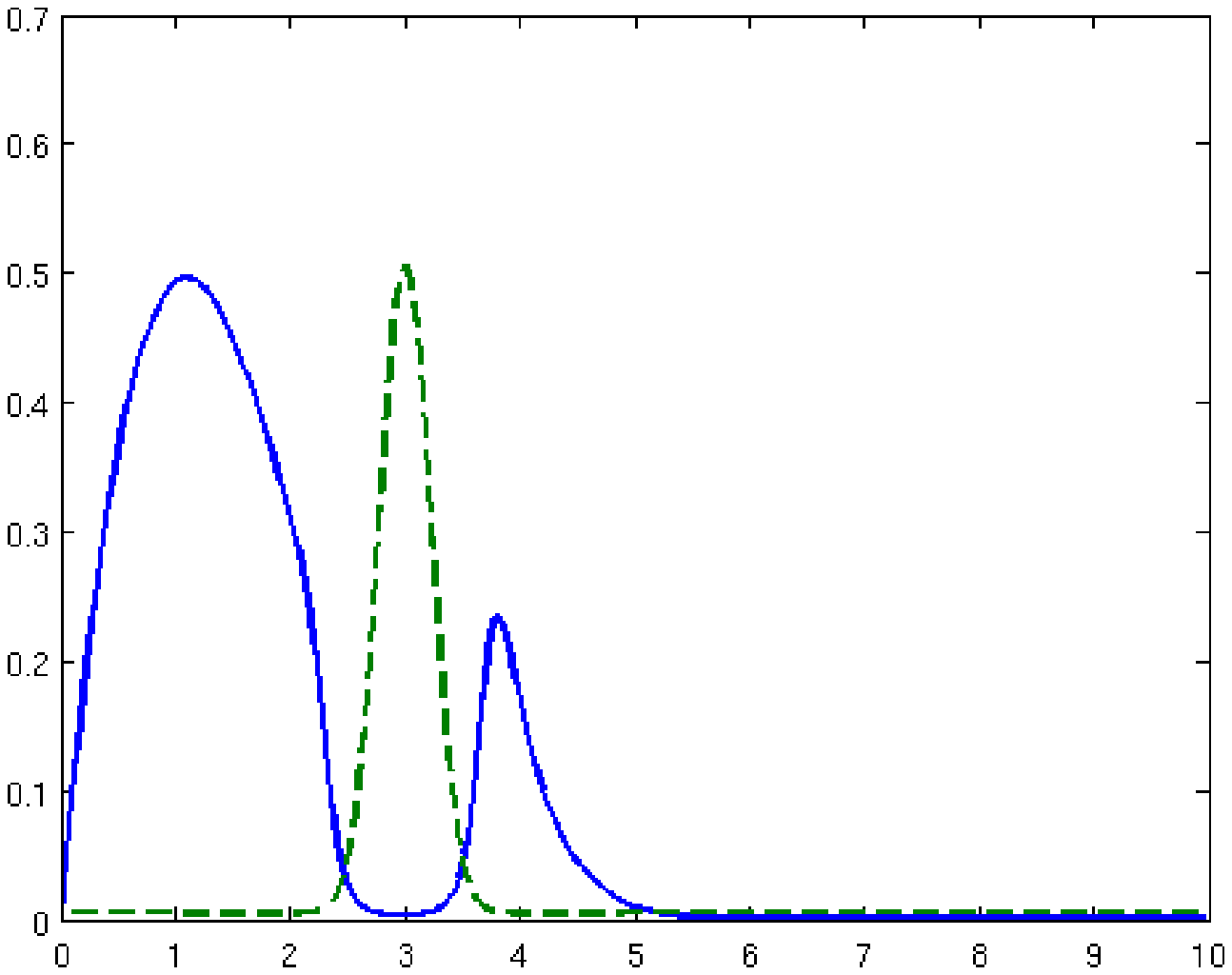}
\includegraphics[width = .45\linewidth]{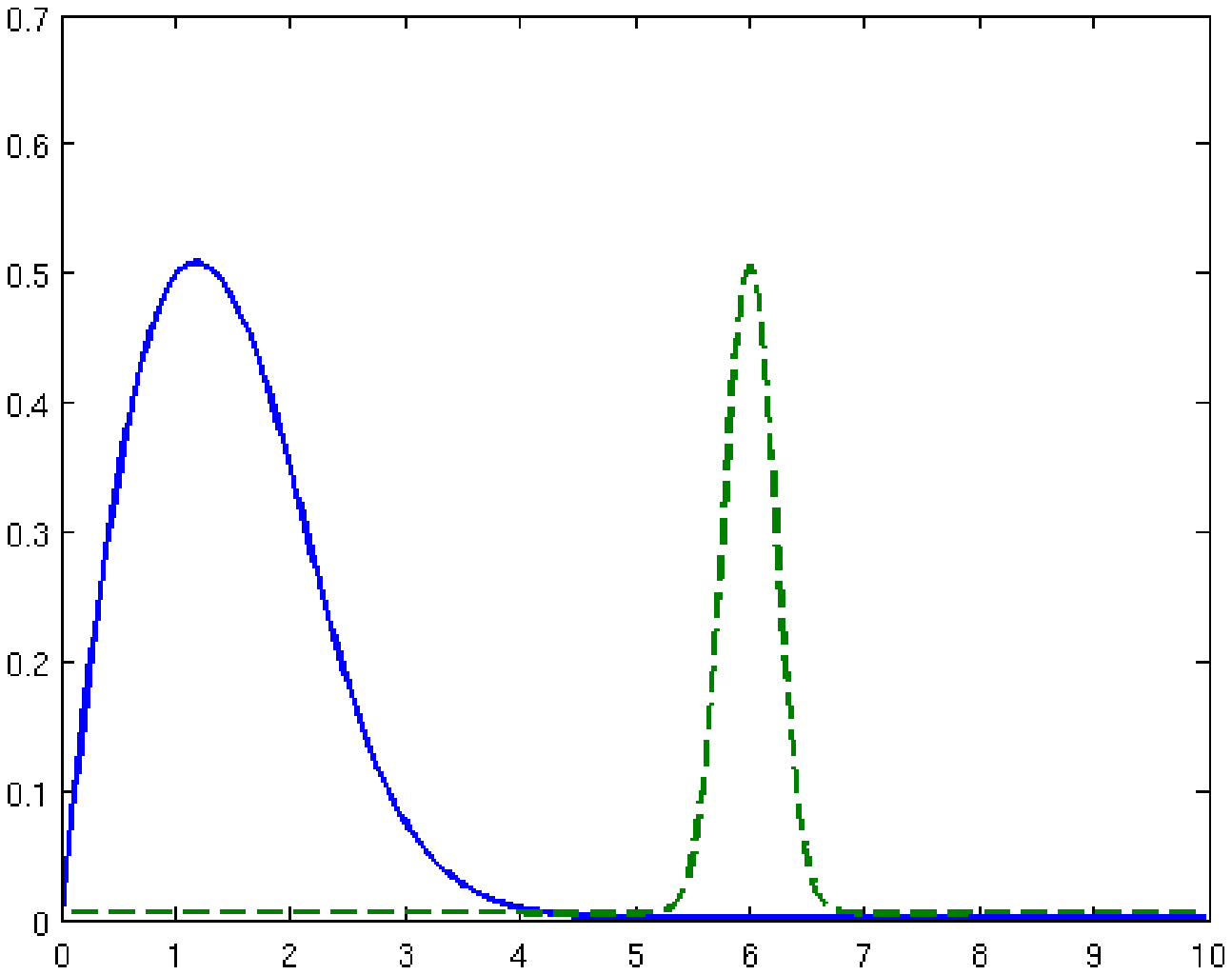}
\caption{{\sc Influence of the peak locus.}
We superpose successively the function $\tau(x)$ (full lines) and the corresponding equilibrium distribution (crossed lines) for translated bells. In each of these pictures, the center of mass of the polymer distribution is 1. We can visualize the splitting effect of the {\em potential} $\tau(x)$ very well.
}
\label{fig:heuristic}
\end{center}
\end{figure}

\paragraph{Locus of the peak of $\tau(x)$: heuristics.} 
For coefficients like $\mu(x) = \mu_0$ and $\beta(x) = \beta_0 x $, we know from  \eqref{eq:Lambda}, that the center of mass at the equilibrium configuration has the value $\mu_0/\beta_0$.  
Departing with $\tau_0$ constant, the unimodal distribution is split into a bimodal one under the action of $\tau(x)$. The further is the peak of the bell-like function $\tau(x)$ from $\mu_0/\beta_0$, the weaker is the splitting effect on the unimodal distribution. This is in accordance with numerical results (see Figure \ref{fig:heuristic}).

\medskip

\section{Stability}
\label{sec:stability}

We have seen previously that there are two possible steady states, setting the alternative between a disease free  and an infected system and depending upon the prion production and degradation rates  $\lambda$, $\gamma$. For 'constant coefficients', Greer {\em et al.} \cite{GreerPW} could study their stability using that the system can be reduced to a differential system on $(V(t), U(t), P(t))$ (see Section \ref{sec:model}). Here, we investigate the same question for general coefficients. Our main assumption is some monotonicity of the eigenvalue in problem  (\ref{eq:frag}),  which we have checked for 'constant coefficients', see formula (\ref{eq:eigenvaccc}). We choose $x_0 = 0$ throughout this section.

\subsection{Stability of the zero state for  $\overline V < V_\infty$}

We first tackle the stability of the disease-free steady state. Recall the notation $\overline V = \frac{\lambda}{\gamma}$ for the 'zero steady state' from (\ref{eq:zerostst}). 

\begin{theorem}[Local stability]  Assume that $\Lambda(\cdot)$ defined in (\ref{eq:frag}) is decreasing 
and that $\overline V < V_\infty$, where $V_\infty$ is defined by $\Lambda(V_\infty)=0$ according to  (\ref{eq:stabless}). Then, in equation (\ref{eq:Greer}), the zero steady state $ V = \overline V $, $u \equiv 0$ is locally nonlinearly stable.
\label{lm:stab1}
\end{theorem}

We recall that in the case at hand ($\overline V < V_\infty$), there does not exist a non-zero steady state because (\ref{as:stst}) cannot be fulfilled with $\varrho >0$. We give in the  subsection \ref{sec:ex} below a class of coefficients for which we can compute the eigenvalue $\Lambda(\cdot)$ and it is indeed decreasing.

\begin{proof} 
The eigenvalue $\Lambda(\cdot)$ is assumed to be a decreasing function of $V$. So, the condition $\overline V < V_\infty$ ensures that $\Lambda(\overline V)>0$. 
We consider a perturbation of the ground state $V(t) = \overline V + \tilde V(t)$ and $u(x,t) = 0+ \tilde u(x,t)$ (note that $\tilde u$ is nonnegative, whereas $\widetilde V$ has no sign).

The  nonlinear system for this perturbation  writes
\begin{equation}
\label{eq:perturbation}
\left\{
\begin{array}{rcl}
\dfrac{d\widetilde V(t)}{dt}&=&\displaystyle - \widetilde V \left[  \gamma + \int_{0}^{\infty} \tau(x) \widetilde u(x,t) \; dx \right] - \overline V \int_{0}^{\infty} \tau(x) \widetilde u(x,t) \; dx ,
\vspace{.2cm}\\
\dfrac{\partial}{\partial t} \widetilde u(x,t) &=& \displaystyle - \overline V \frac{\partial}{\partial x} \big(\tau(x) \widetilde u(x,t)\big)- \widetilde V(t) \frac{\partial}{\partial x} \big(\tau(x) \widetilde u(x,t)\big) - [\mu(x) + \beta (x)] \widetilde u(x,t)
\vspace{.2cm}\\
&&  \qquad + \displaystyle2 \int_{x}^{\infty}\, \beta(y) \kappa(x,y)  \widetilde u(y,t) \, dy. 
\end{array} \right.
\end{equation}

Following the duality  method in \cite{CCP}, we introduce the adjoint eigenvector for the eigenvalue $\overline V$, namely $\overline \varphi >0$ given by the equation 
\begin{equation}
\label{eq:adjoint}
\displaystyle - \overline V  \tau(x)  \frac{\partial}{\partial x} \overline \varphi(x) +[ \mu(x) +  \beta (x)]\overline \varphi(x)  = 2 \int_{0}^{x}\beta(x) \kappa(y,x) \overline \varphi(y) \, dy + \Lambda(\overline V) \overline \varphi(x)\ . 
\end{equation}
We now introduce 
constants $K_1,\, K_2$ such that
\begin{equation}
\left| \tau(x) \;\dfrac{\partial}{\partial x}\overline \varphi(x)  \right| \leq K_1 \overline \varphi(x) , \qquad \tau(x) \leq K_2 \overline \varphi(x)\, .
\label{as:stab}
\end{equation}
This is possible because $\overline  \varphi$ grows linearly at infinity according to general abstract properties proved in  \cite{ Perthame_LN,michel1, PeRy}. See again subsection \ref{sec:ex} for examples.  

We test the adjoint equation against the equation on $\widetilde u$ in \eqref{eq:perturbation} and obtain
\[ \dfrac d {dt} \int_{0}^{\infty} \widetilde u(x,t) \overline \varphi(x)\, dx = - \Lambda(\overline V)\int_{0}^{\infty} \widetilde u(x,t) \overline \varphi(x)\, dx + \widetilde V(t)\int_{0}^{\infty} \left( \frac{\partial}{\partial x} \overline \varphi(x) \right) \tau(x) \widetilde u(x,t)\, dx \ . 
\]

On the other hand, multiplying the first differential equation in \eqref{eq:perturbation} by the sign of $\widetilde V$, we get:
\[ \dfrac {d}{dt}|\widetilde V| \leq 
- |\widetilde V| \left[  \gamma + \int_{0}^{\infty} \tau(x) \widetilde u(x,t) \; dx \right]+\overline V \int_{0}^{\infty} \tau(x) \widetilde u(x,t) \; dx\, . \]

We obtain, choosing $\alpha$ large enough such that  $\delta \Lambda(\overline V) -\frac {K_2 \overline V}{ \alpha} >0 $, 
\begin{eqnarray*}  
\dfrac d{dt}\left( \alpha  \int \widetilde u \overline \varphi + |\widetilde V| \right) 
 &\leq&  
- \alpha \Lambda(\overline V)\int \widetilde u \overline \varphi 
 + \alpha K_1  |\widetilde V(t)|  \int  \overline \varphi \widetilde u -  \gamma |\widetilde V|  - |\widetilde V| \int  \tau  \widetilde u + K_2 \overline V \int \widetilde u \overline \varphi  \\
&\leq&  
- \min\left(\Lambda(\overline V) -\frac {K_2 \overline V}{ \alpha} , \gamma\right) \left(\alpha  \int \widetilde u \overline \varphi 
+  |\widetilde V|\right) + \alpha K_1  |\widetilde V(t)|  \int  \overline \varphi \widetilde u 
 \\
&\leq&  -\min(\delta,\gamma) \left(\alpha  \int \widetilde u \overline \varphi 
+  |\widetilde V|\right) + \frac{K_1}{2} \left(\alpha  \int \widetilde u \overline \varphi 
+  |\widetilde V|\right)^2\ .
\end{eqnarray*}
From this differential equation we conclude that, when $\left(\alpha  \int \widetilde u \overline \varphi 
+  |\widetilde V|\right)$ is initially small enough, then the right hand side is negative. Therefore it decays for all times with the asymptotic exponential rate 
$$ \left(\alpha  \int \widetilde u \overline \varphi +  |\widetilde V|\right) \lesssim C e^{-\min(\delta,\gamma) t} \ .
$$
\end{proof}

We can also state the following variant of Theorem \ref{lm:stab1}

\begin{theorem}[Global stability]  Additionally to the hypotheses of Theorem \ref{lm:stab1}, assume   that  $\tau(x) \geq k \overline \varphi(x)$ for some constant $k>0$ and $\overline V/\Lambda(\overline V)$ is small enough compared to $k/{(K_1K_2)t}$. Then, in equation (\ref{eq:Greer}), the zero steady state $ \overline V= \frac{\lambda}{\gamma}$, $u \equiv 0$ is globally nonlinearly stable.
\label{lm:stab2}
\end{theorem}

By opposition to those of Theorem \ref{lm:stab1}, the assumptions of this Theorem are difficult to check directly because the coefficients are intricated here. They mean that we are close to the 'constant coefficients case' because $\overline \varphi = 1$ in this case and we can choose $K_1 = 0$.

\begin{proof} With these additional assumptions, we may keep one negative term in the last computation and arrive to 
\begin{eqnarray*}  
\dfrac d{dt}\left( \alpha  \int \widetilde u \overline \varphi + |\widetilde V| \right) 
 &\leq&  
- \alpha \Lambda(\overline V)\int \widetilde u \overline \varphi 
 + \alpha K_1  |\widetilde V(t)|  \int  \overline \varphi \widetilde u -  \gamma |\widetilde V| - |\widetilde V|\int \tau\widetilde u  + K_2 \overline V \int \widetilde u \overline \varphi  \\
&\leq&  
- \min\left(\Lambda(\overline V) -\frac {K_2 \overline V}{ \alpha} , \gamma\right) \left(\alpha  \int \widetilde u \overline \varphi 
+  |\widetilde V|\right) +( \alpha K_1- k)  |\widetilde V(t)|  \int  \overline \varphi \widetilde u
 \\
&\leq&  -\min\left(\Lambda(\overline V) -\frac {K_2 \overline V}{ \alpha} , \gamma\right) \left(\alpha  \int \widetilde u \overline \varphi +  |\widetilde V|\right).
\end{eqnarray*}

This last inequality will hold if we can find $\alpha > \frac{K_2 \overline V}{\Lambda(\overline V)}$ such that $\alpha K_1 <k$, which is precisely our smallness assumption.  Then we have exponential decay of the solution.
\end{proof}

\subsection{Unstability of the zero state for  $\overline V > V_\infty$}

The same kind of method allows for another type of results initiated in \cite{CCP}, namely that solutions cannot go extinct when another steady state exists. Again we choose $x_0=0$ in this subsection. We have, recalling that $\overline  \varphi$ is defined in (\ref{eq:adjoint}),  the 

\begin{theorem}[Unstability]  Assume that $\Lambda(\cdot)$ defined in (\ref{eq:frag}) is decreasing, that for some constant $K>0$, $ \tau(x) \;\partial_x\overline \varphi(x) \leq K \overline \varphi(x)$, that  $V(0) \leq \overline V$  and that $ V_\infty < \overline V$. Then, in equation (\ref{eq:Greer}), the zero steady state $ V = \overline V $, $u \equiv 0$ is locally unstable.
\label{lm:unstab1}
\end{theorem}

\begin{proof} Following the proof of Theorem \ref{lm:stab1}, we have on the one hand
$$
\frac {d}{d t} (\overline V -V(t)) + \gamma (\overline V -V(t))=  V(t) \int_0^\infty \tau(x) u(x,t)\, dx \geq 0.
$$
Therefore we have $V(t) \leq \overline V$ for all $t\geq 0$. 

On the other hand, we can again combine the equations (\ref{eq:adjoint}) for   $\overline \varphi$  and the equation on $u$    in (\ref{eq:Greer}), to obtain
$$
\frac {d}{d t} \int_0^\infty u(x,t) \overline \varphi(x)\, dx  = (V(t)-\overline V) \int_0^\infty \tau(x) \dfrac{\partial}{\partial x} \overline \varphi(x) \; u(x,t)\, dx -\Lambda(\overline V) \int_0^\infty  u(x,t) \overline \varphi(x)\, dx.
$$

We define the quantity (sum of two positive terms)
$$
Q(t)= \int_0^\infty u(x,t) \overline \varphi(x)\, dx + \overline V -V(t).
$$
As a consequence of the above calculations and of the hypothesis $ \tau \;\partial_x\overline \varphi  \leq K \overline \varphi $, we have therefore,
\begin{eqnarray*}
\frac {d}{d t} Q(t)& \geq& - K (\overline V - V(t)) \int_0^\infty \overline \varphi(x) \; u(x,t)\, dx - \Lambda(\overline V) \int_0^\infty   u(x,t) \overline \varphi(x)\, dx
 \\ 
&\geq&  k_1 [1- k_2 (\overline V - V(t))] \int_0^\infty   u(x,t) \overline \varphi(x)\, dx \,  ,
\end{eqnarray*}
for some constants $k_i>0$. Indeed, in this case when $V_\infty \leq \overline V$, we know that $ \Lambda(\overline V)<0$.

We can  conclude from this inequality that $0$ is unstable. Indeed, whenever $Q(t)$ becomes small enough, then 
$\overline V - V$ becomes also small enough so that $Q(t)$ increases and the solution gets away from $0$. 
\end{proof}
 
\subsection{A class of examples}
\label{sec:ex}

In order to clarify the assumptions and  properties stated before, we give a class of coefficients in equation (\ref{eq:Greer}) where the different eigenelements can be explicitly computed. In each case we will see that $
\Lambda(V)$ s a decreasing function of $V$, and that (\ref{as:stab}) reduces to the fact that $\tau(x)$ is bounded. 
\\

We first recall from \cite{Perthame_LN, PeRy}, the case $\tau\equiv\tau_0$, $\mu\equiv\mu_0$ and $\beta\equiv\beta_1$ constant. We obtain the solution to (\ref{eq:adjointp}) below,
$$
\Lambda(V)= \beta_1, \qquad \varphi(x)=\text{constant}.
$$

The eigenlements are usually difficult to evaluate in the direct problem (\ref{eq:frag}) but can be easily computed in the adjoint equation, recall  (\ref{eq:adjoint}), 
\begin{equation}
\label{eq:adjointp}
\displaystyle - V\,   \tau(x)  \frac{\partial}{\partial x}  \varphi(x) +[ \mu(x) +  \beta (x)] \varphi(x)  = 2 \int_{0}^{x}\beta(x) \kappa(y,x)  \varphi(y) \, dy + \Lambda( V)  \varphi(x)\ . 
\end{equation}
Indeed, searching for an affine  solution $\varphi(x) =1+x/L$, we find, using the structure properties of $\kappa(x,y)$ in (\ref{as:kappa}),
$$
0= \frac{V \tau(x)}{L}+ \beta(x)  +(\Lambda-\mu(x))\left(1+\frac x L\right).
$$

Hence a more general class where one can compute $\Lambda(V)$ is when $\beta(x)= \beta_0\, x$, and $\tau_0$, $\mu_0$ are constant. We arrive at
$$
-\beta_0 L  ^2  + \tau_0 \,  V=0.
$$
Therefore we obtain $L = \sqrt{\tau_0/\beta_0}$ and 
$
\Lambda(V)= \mu_0  - \beta_0  L( V),
$
 as in Section \ref{sec:model}.
Observe that $\Lambda(V)$ is a decreasing function of $V$ as used in Theorem \ref{lm:stab1}.

We leave to the reader to check that $\Lambda(V)$ is also decreasing when $\beta(x)=\beta_1 + \beta_0\, x$, $\tau(x)= \tau_0 + \tau_1\, x$ and $\mu_0$ is constant. Indeed solving the two equations for $\varphi(x) =1+x/L$, we find 
$$\Lambda-\mu_0  +\frac{ V\tau_0}{L} + \beta_1 =0\, , \qquad \frac{\Lambda-\mu_0 } L  + \frac{V \tau_1}L+  \beta_0 =0. 
$$
And thus $Z=\Lambda-\mu_0 $ solves
$$
(Z +\beta_1)( Z + V\tau_1)=  V \tau_0   \beta_0\, , \qquad Z < - V\tau_1.
$$
This last constraint comes from the condition $L>0$.
\\

More generally, notice that it always holds true that 
\begin{equation} \label{eq:frag1}
\Lambda(0) =  \mu_0\, , \qquad \Lambda(V) \le  \mu_0\, .  
\end{equation}
To prove it, we just multiply  (\ref{eq:frag})  by $x$  and integrate. Notice that the problem is then degenerate and the eigenfunction is singular at $x=0$.

On the other hand, still for the same reason, if $\min \frac{\tau(x)}{x} >0$, we have 
\begin{equation} \label{eq:frag2}
\Lambda(\infty) =  - \infty, 
\end{equation}
and then there is at least  one  non-zero steady state.

\section{Temporal dynamics of the aggregates size distribution in a biological context}
\label{sec:timescales}
%

As we highlighted previously, recent experimental work by Weber et al \cite{Weber1, Weber2} and by Silveira et al \cite{Caughey} suggests that PrPsc aggregates are differentially infectious with respect to their size. The former also shows that the PrPsc size distribution in a whole infected brain is bimodal rather than unimodal. Based on the nucleated polymerization model for prion growth, we have shown that bimodal stationary distribution can occur for a size-dependent conversion rate. However, during the time course of incubation period and clinical stage of prion diseases, PrPsc accumulation seems to follow an exponential dynamic until the death. Apart from hemizygote mice (PrP+ /0), no steady state of PrPsc accumulation was observed in the brain during all the disease \cite{Bueler, Collinge}.  
Consequently, we focus in this section on the dynamic of PrPsc accumulation in the early stages of the infection named 'exponential phase'. 
Our approach is based both on the eigenvalue problem and on direct numerical simulations of the temporal dynamics.

\paragraph{Numerics.}
Parameters of the nucleated polymerization model for prion growth have been estimated for the 'constant parameters model' from experimental data and are available in the literature \cite{masel, GreerPW, Rubenstein}. The parameter values used in the sequel have been   quoted from Rubenstein et al. \cite{Rubenstein}. Unless explicitly mentioned, they are:
$ \lambda = 2400 $ per day, $ \gamma= 4$ per day, $\mu_0=.05 $ per day, and $ \beta_0= 0.03$ per day. The conversion function $\tau(x)$ is the sum of a basal rate $\tau_0 = .001$ and a gaussian bell centered on $m$:
\begin{equation} \tau(x)=\tau_0+A\exp\left(-(x-m)^2/\sigma^2\right)\, , 
 \label{eq:bell shape}
\end{equation}
 with a magnitude $A$ to be chosen by several orders of magnitude above $\tau_0$ (from biological evidence on different specific converting activities in \cite{Caughey}).
The simulations assume an initial PrPc population $V(0) = \overline V = \lambda / \gamma $ (corresponding to the healthy steady state) and an initial PrPsc population given by $u(x,0) = 0.5  x^2 /(1+x^4) $, which is a small perturbation of the zero steady state indeed.

\paragraph{Dynamics of PrPsc accumulation and Eigenvalue problem.}
The eigenvalue problem can be used to measure the dynamics of PrPsc aggregates accumulation during the exponential growth.
We assume below that $\overline V>V_\infty$, that $\Lambda(V)$ is a decreasing function and that instability of the healthy state holds true (see Section \ref{sec:stability}). The practical choice of coefficients is mentioned above.

Assuming that $V(t)$ remains close enough to $\overline V$ (this is exponential growth phase in a linear regime), then the polymerization behaves following a linear problem. Actually the second equation in \eqref{eq:Greer} is decoupled from the first one at the first order of approximation. The dominant eigenvalue $-\Lambda(\overline V)$ thus measures the exponential growth of the PrPsc total population (see \cite{Perthame_LN} for a mathematical formulation of this fact using the generalized relative entropy). As a by-product, we also learn that the renormalized distribution tends to align along the eigenvector $\mathcal U(\overline V,x)$. One noticeable difference with the stationary case studied in Section \ref{sec:model} is that the center of mass of $\mathcal U(\overline V,x)$ is translated to the right as opposed to $\mathcal U(V_\infty,x)$ (for which it is $\mu_0/\beta_0$), according to \eqref{eq:Lambda}:
\[0< - \dfrac{\Lambda(\overline V)}{\beta_0} = \int \left( x - \dfrac{\mu_0}{\beta_0} \right) \mathcal U(\overline V,x)\, dx\, .  \] Observe that for the parameters picked up in the literature, the difference is significant (approximately 500\%).

\medskip

\begin{figure}
\begin{center}
\includegraphics[width = .6\linewidth]{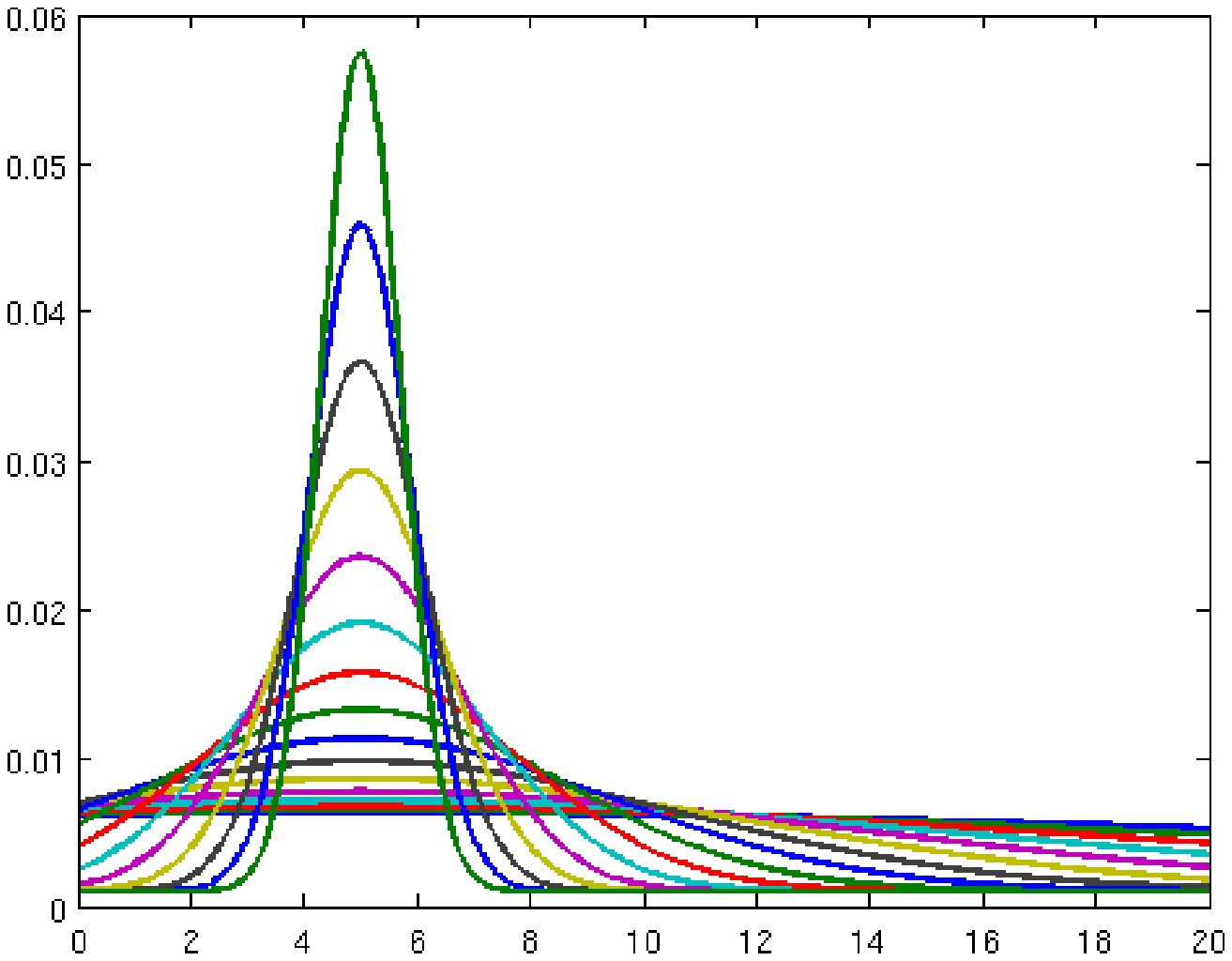} \\
\includegraphics[width = .6\linewidth]{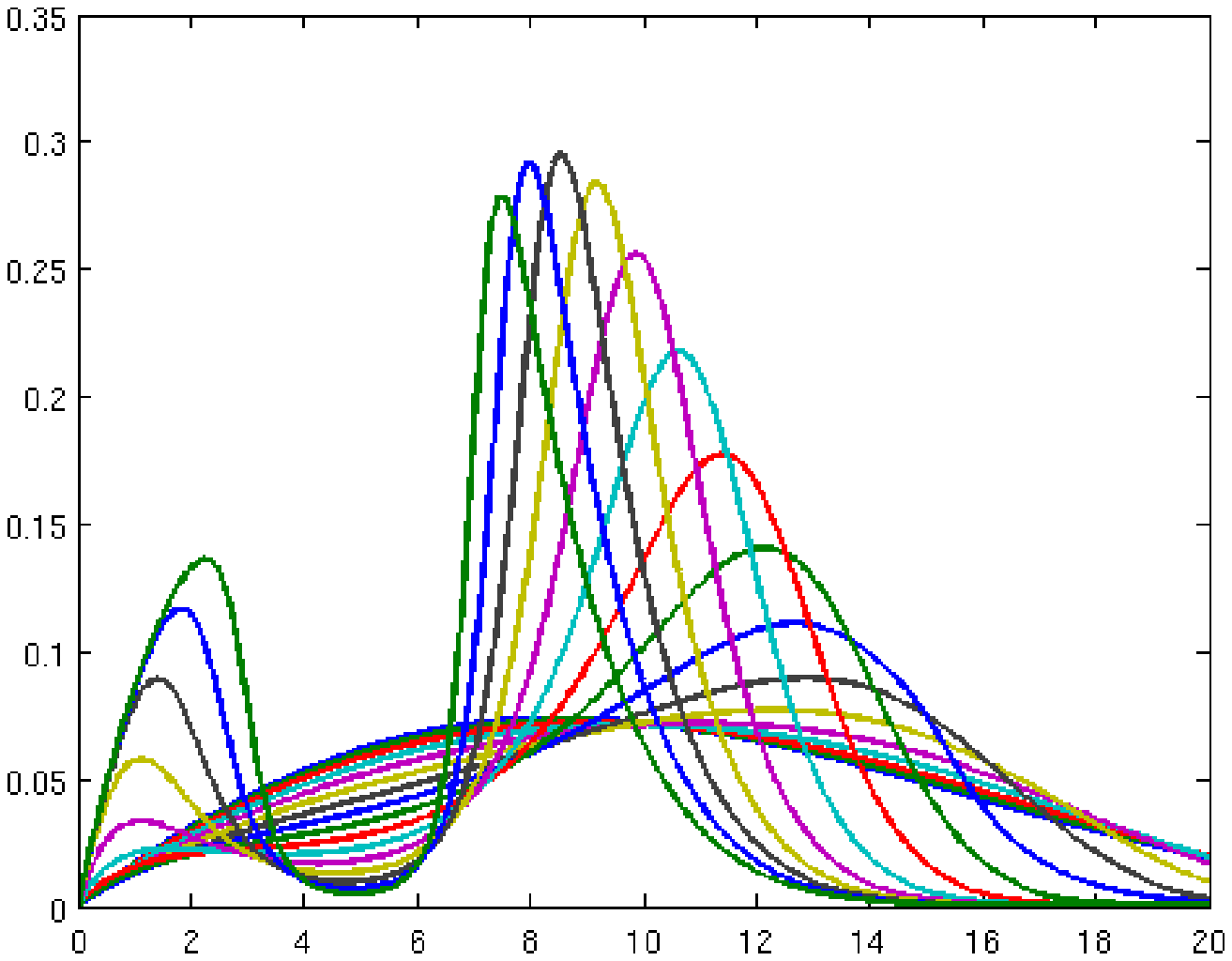} \\
\includegraphics[width = .6\linewidth]{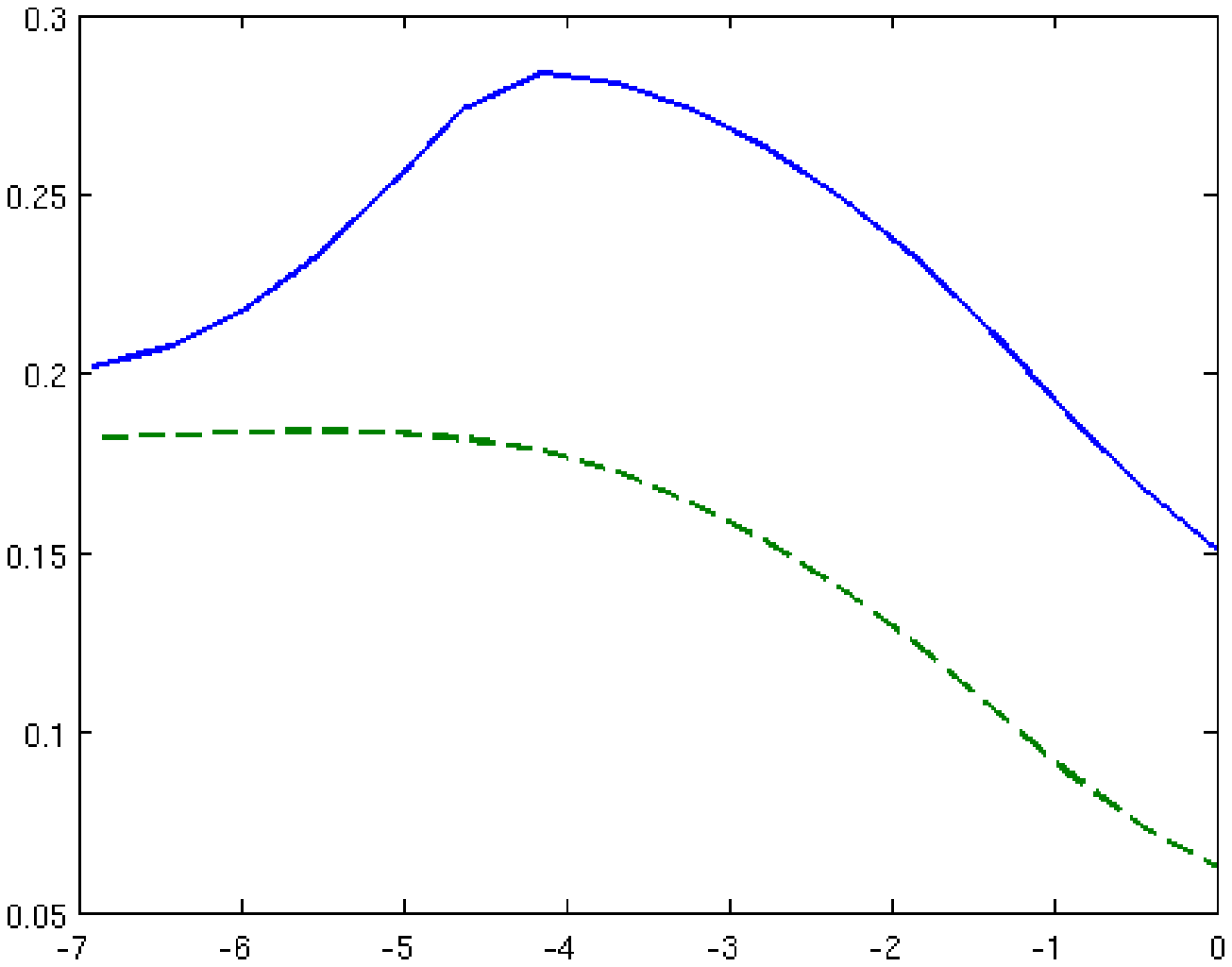} 
\caption{{\sc Influence of the transconformation tightness} (top) The transconformation rate with several levels of concentration: $\tau(x) = \tau_0 + \alpha\varphi(\alpha(x-m))$, where $\varphi$ is a gaussian function and $\alpha = 10^{(-3:.2:0)}$). (middle) The corresponding eigenfunctions $\mathcal U(\overline V;x)$. (bottom) The exponential growth rate $-\Lambda(\overline V)$ (solid line) as a function of $\alpha$ (logarithmic scale, units for $-\Lambda(\overline V)$ are day$^{-1}$) and the effective transconformation rate (dashed line) $\tau_{\mbox{\small eff}} = \int \tau(x)\mathcal U(\overline V,x)\, dx$ (magnified 30 times). Bimodal distribution begins for $\ln(\alpha) \gtrsim -3$ and thus does not correspond to the fittest distribution.}
\label{fig:vp}
\end{center}
\end{figure} 

\medskip

If $\varrho(t)$ denotes the total population of PrPsc, then it evolves eventually as $\varrho(t)\approx \varrho(0) \exp\big(-\Lambda(\overline V) t \big)$ from the ground healthy state \cite{Perthame_LN} in the exponential phase. Here the initial condition $\varrho(0) = \varrho_{\mbox{\small inoculation}}$ accounts for the inoculation dose. Therefore $\log \varrho(t)$ grows linearly in time as a first approximation. If we define mathematically the incubation time as a threshold in that macroscopic amount of PrPsc \cite{GreerPW}, we rederive the logarithmic relationship between incubation time and inoculation dose in the linear regime:
\begin{equation} T_{\mbox{\small incubation}} 
\approx -\dfrac1{\Lambda(\overline V)}  \log\left(\dfrac{\varrho_{\mbox{\small symptomatic}}}{\varrho_{\mbox{\small inoculation}}}\right) \, . \label{eq:log law} \end{equation}
Moreover, the relationship \eqref{eq:log law} stresses out the importance of measuring the dependence of $\Lambda(\overline V)$ with respect to the parameters of the model. Following Figure \ref{fig:eigenvalue} we suggest to evaluate the influence of the tightness of the transconformation peak. Closely looking at Figure \ref{fig:eigenvalue}(right) we indeed notice that the growth rate $-\Lambda(\overline V)$ (right limit of the graph) is slightly larger for the model with constant rate $\tau(x) = \tau_0$ (solid line). Several transconformation rates $\tau(x)$ -- being the same basal rate $\tau_0$ combined with a more and more concentrated gaussian bell are tested -- and the corresponding growth rate in the exponential phase $-\Lambda(\overline V)$ are computed numerically (Figure \ref{fig:vp}). Interestingly, the results exhibit a best compromise around $\alpha \approx 0.01$ (intermediate concentrations of the peak). However, this does not correspond to a bimodal distribution for the polymerization profile. Thus, according to this model the optimal conditions for PrPsc accumulation in the exponential expansion phase are not those which lead to a bimodal size-distribution of PrPsc aggregates as observed by Silveira et al. \cite{Caughey}.

\paragraph{Prion strain mechanisms.}

\begin{figure}
\centering
\includegraphics[width=.3\linewidth, height=.3\linewidth]{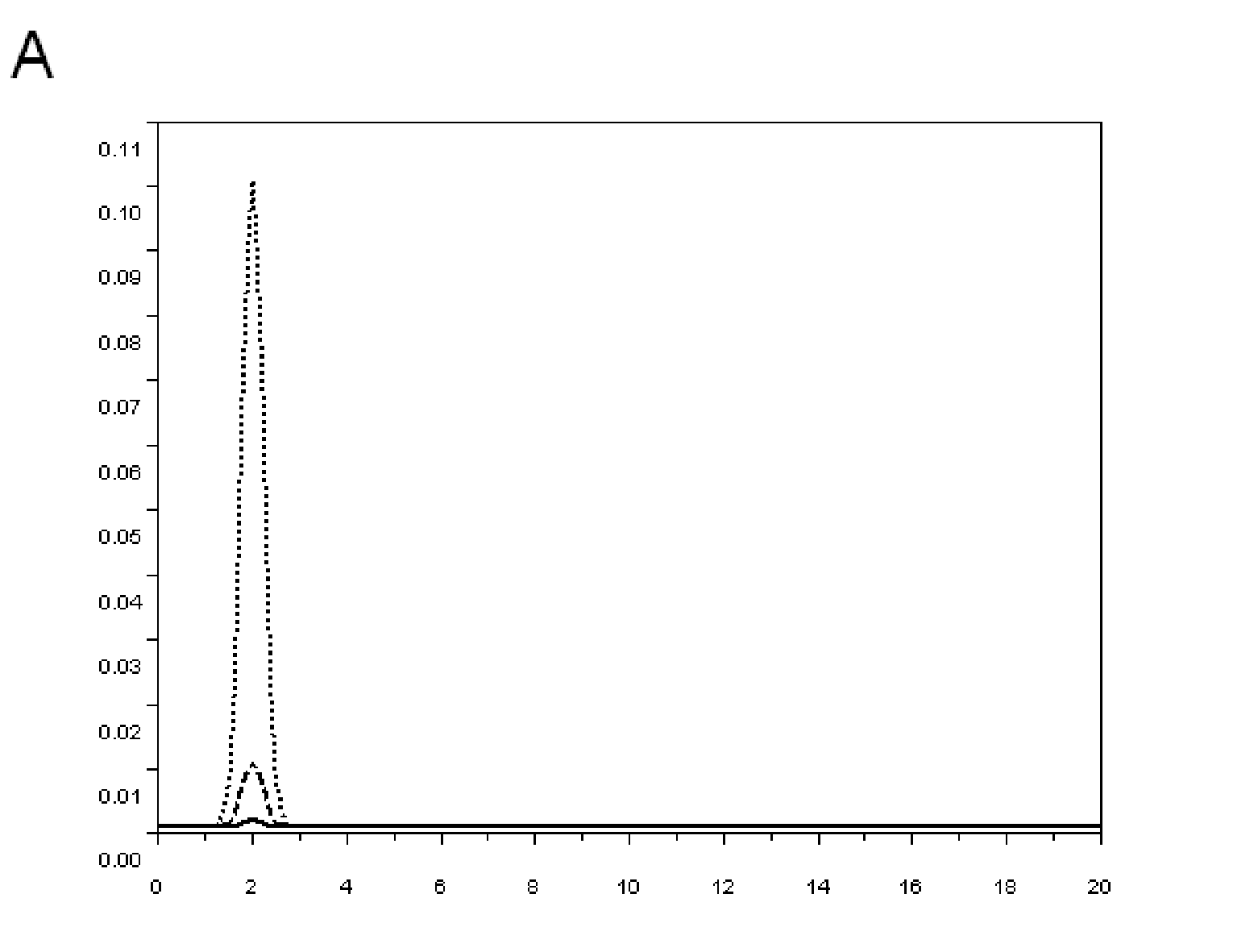} 
\includegraphics[width=.3\linewidth, height=.3\linewidth]{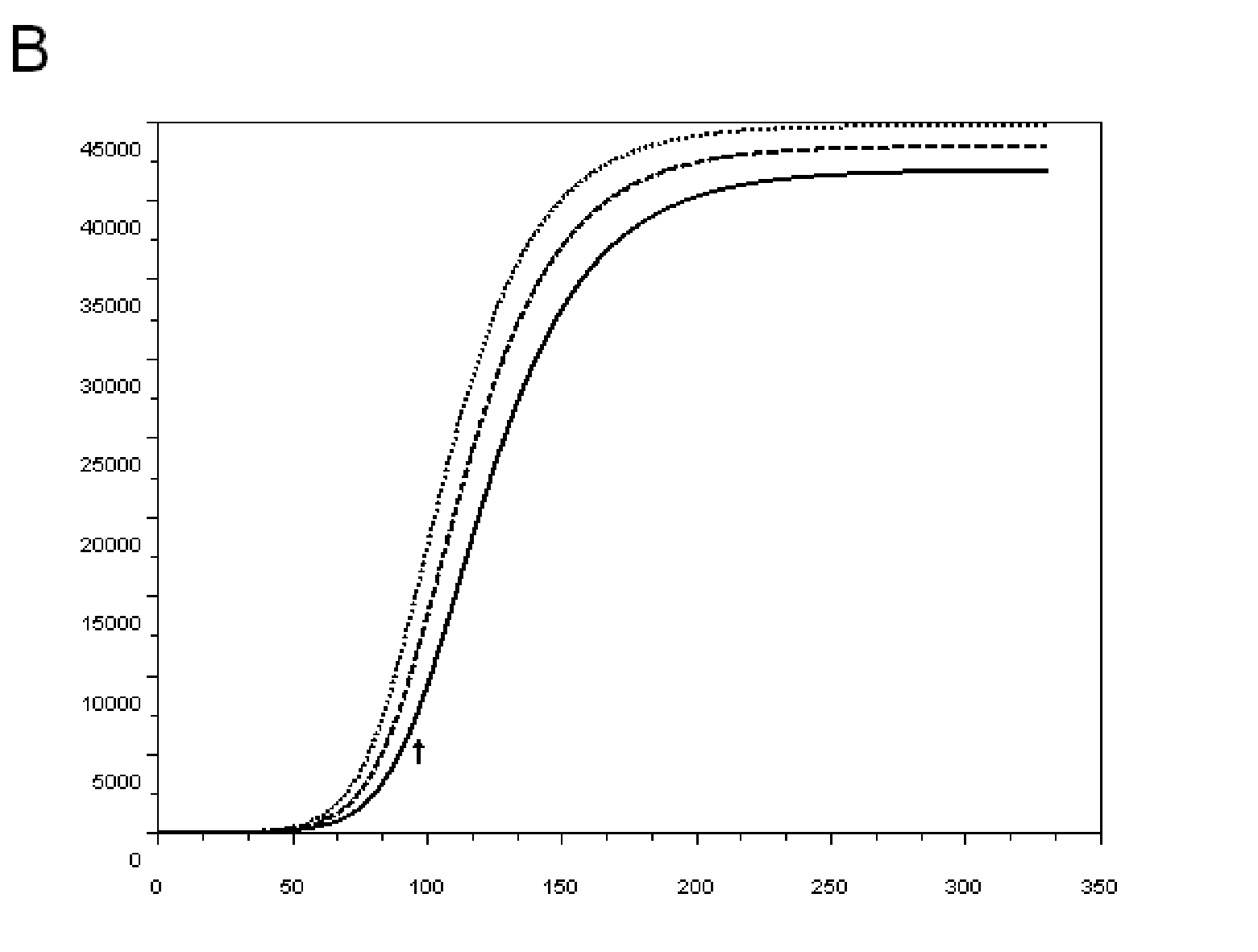} 
\includegraphics[width=.3\linewidth, height=.3\linewidth]{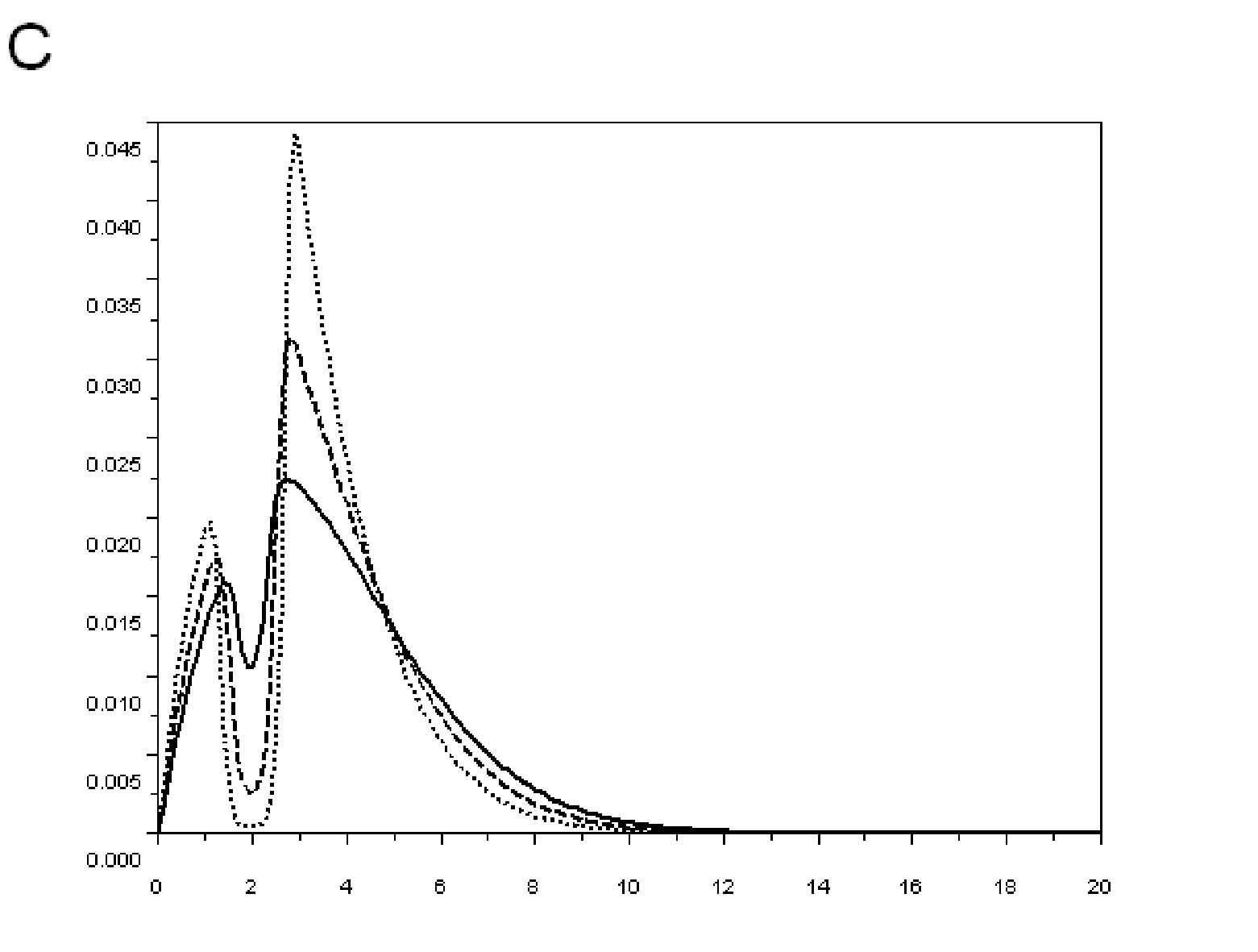} 
\caption{Prion replication with variations in the maximum level of the transconformation rate $\tau$. Transconformation rate is given by $\tau(x) = 0.001 + H*exp(-10*(x-2)^2)$. Three values for $H$ have been tested :  $H = 0.001$ (solid line), $H=.01$ (dashed line) and $H=0.1$ (dotted line). 
(A) Size distribution of $\tau$ (abscissa = PrPsc aggregates size; ordinate =  rate $\tau$)
(B) Time evolution of total PrPsc. The arrow represents time t = 96 days (abscissa = Time (in day) ; ordinate =  rate $\tau$ (per day)
(C) Normalized PrPsc aggregates size distribution at t = 96 days, corresponding to the exponential growth of PrPsc. The distributions are normalized by the total number of PrPsc aggregates (abscissa = PrPsc aggregates size; ordinate =  PrPsc aggregates number). 
} 
\label{fig:htau} 
\end{figure}

\begin{figure}
\centering
\includegraphics[width=.3\linewidth, height=.3\linewidth]{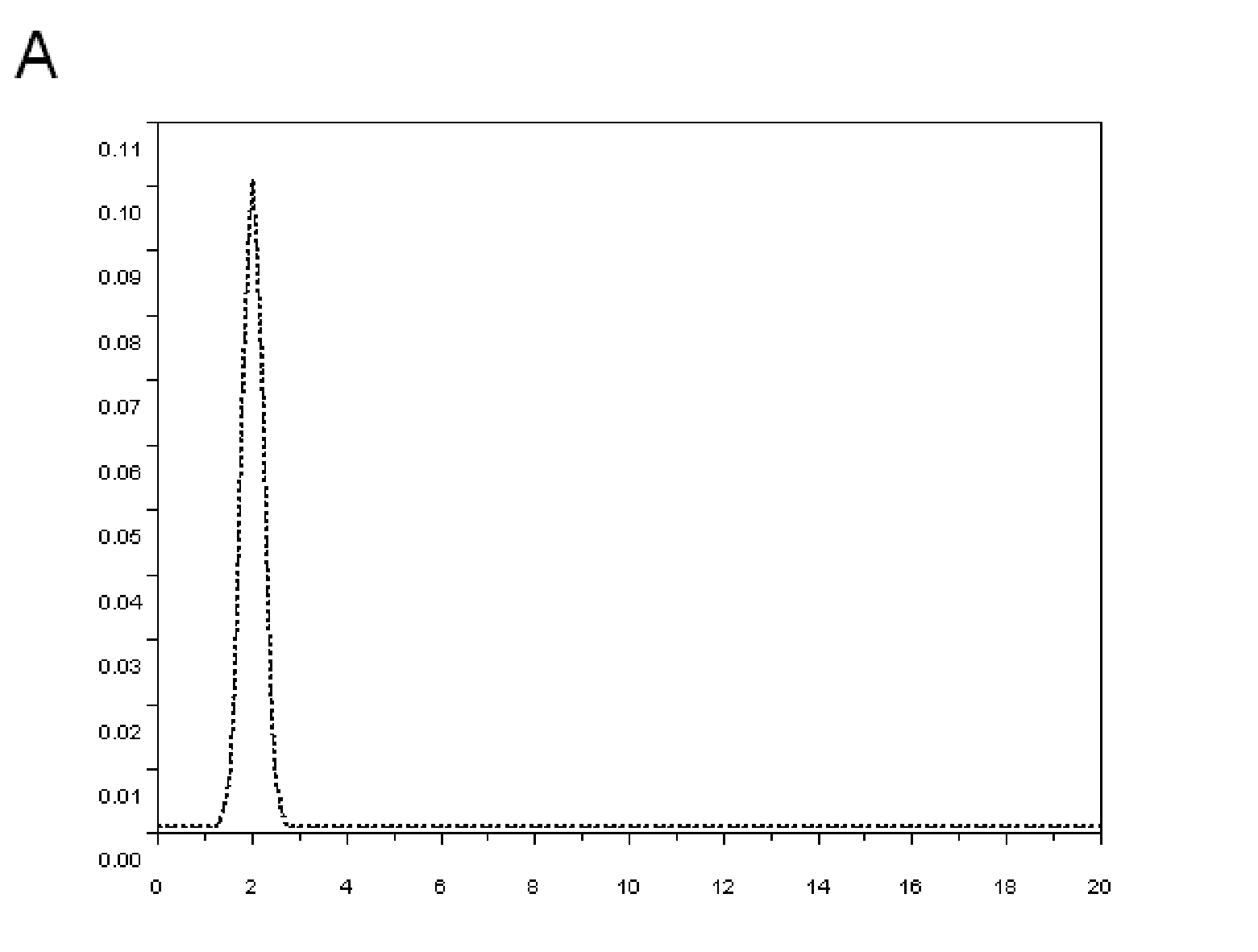} 
\includegraphics[width=.3\linewidth, height=.3\linewidth]{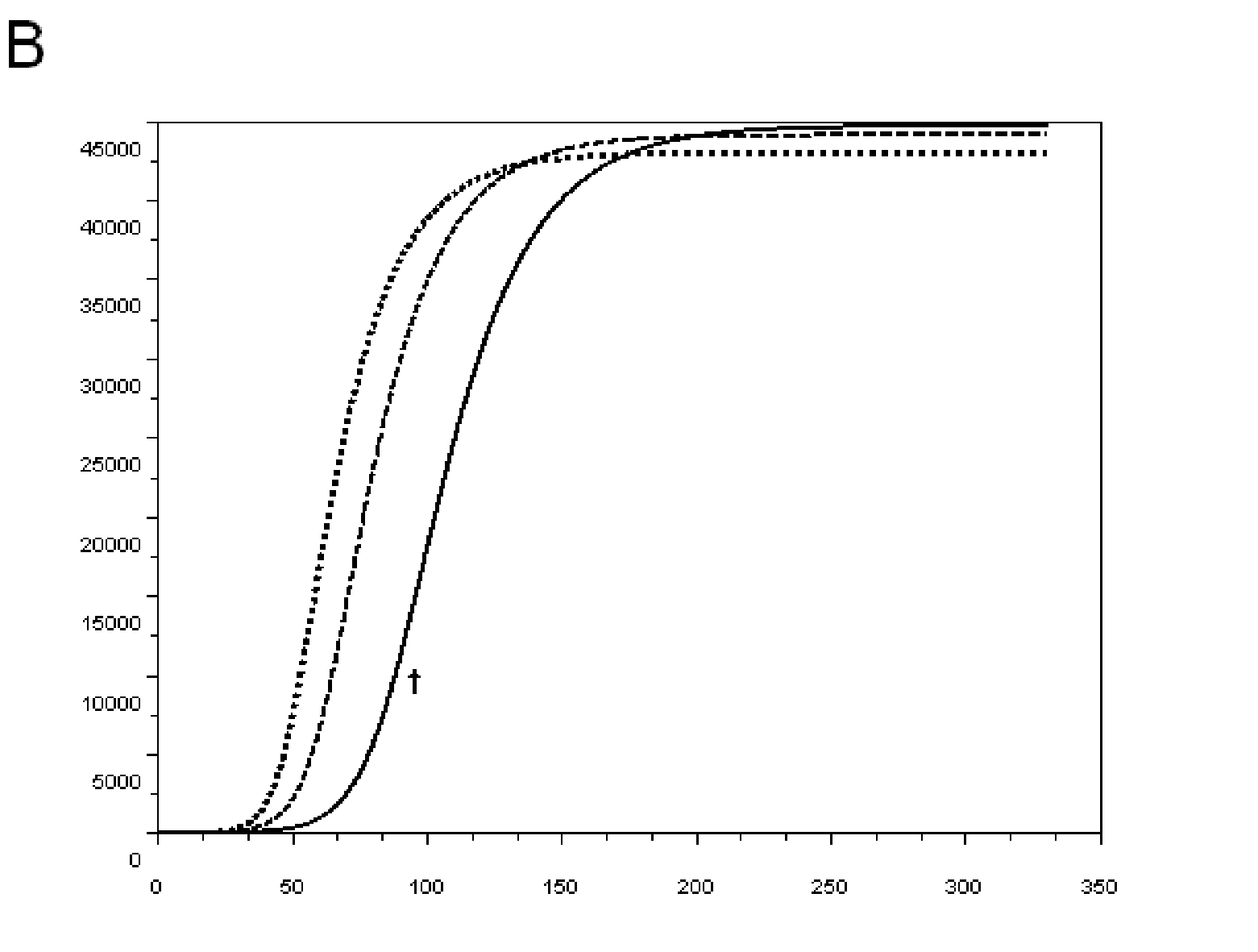} 
\includegraphics[width=.3\linewidth, height=.3\linewidth]{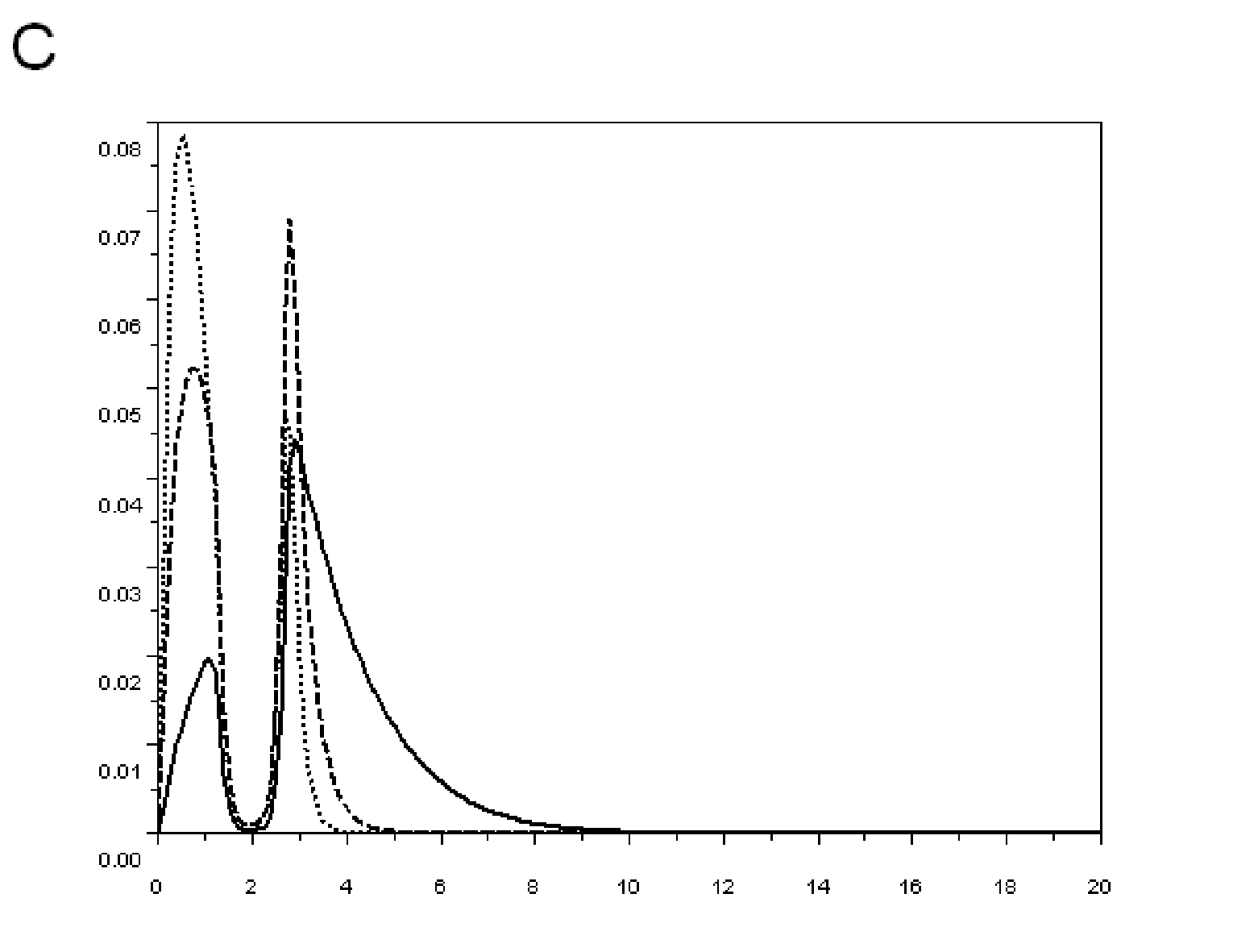}  
\caption{ Prion replication with variations in the fragmentation rate of aggregates. Three value for $\beta$ have been tested :  $\beta = 0.0314$ (solid line), $\beta = 0.0471$ (dashed line) and $\beta = 0.0628$ (dotted line). 
(A), (B) and (C) : Same as Figure \ref{fig:htau}} 
\label{fig:betavar} \end{figure}

It is a challenging problem to investigate the prion strain phenomenon under the light of model \eqref{eq:Greer} which provides informations on the possible microscopic distributions of PrPsc polymers. The prion strain phenomenon is supported by the fact that prion infected animals may develop several distinct pathologies, whose clinical and neuropathological outcomes can be maintained through several passages in rodents models of prion diseases. Each prion strain transmitted in the same congenic host exhibits stable incubation period, clinical expression, lesions profile and infectious properties. These phenotypic features can be used to differentiate between different prion strains \cite{Morales}. Among them, the most commonly used is the incubation period. For most prion strains within the same host, incubation period is determined by the rate and pattern of PrPsc accumulation. 
For our purpose however it will be reduced to the time for the total amount of PrPsc to reach a certain threshold. 
We therefore study by mean of direct numerical simulations how parameters influence the dynamics of PrPsc accumulation. 
The values used for these numerical tests are mentioned above. The locus and the thightness of the transconformation rate are  qualitatively compatible with a bimodal distribution of PrPsc polymers in the exponential growth phase (see figure \ref{fig:htau} and figure \ref{fig:betavar}, B), as well as in the equilibrium phase (data not shown). For the experiments depicted in figure \ref{fig:htau} and figure \ref{fig:betavar}, varying parameters acts upon the kinetics of PrPsc accumulation as follows. The higher the transconformation or the fragmentation rates are, the faster PrPsc accumulates, and thus, the shorter the incubation time seems to be. In particular this implies that a fast strain could be either instable (i.e. can be easily broken) or able to fix and transconform PrPc with great efficiency (or both). Interestingly, our simulations also predict changes in the microscopic size-distribution of PrPsc aggregates, whose profiles are specific of the given parameters.


\section{Conclusion and perspectives}
Our motivation for this work is to better understand the unexpected bimodal size distribution of prions, observed by Silveira et al \cite{Caughey}.
Based on the nucleated polymerization model of prion fibril growth describing the dynamics of PrPc monomeres and PrPsc aggregates, we take into account aggregate-size dependent parameters, expanding thus a previous study \cite{GreerPW}. 
We have shown numerically that bimodal stationary distribution can occur in the case of a non-uniform conversion rate from PrPc into PrPsc. It is worth noticing that experimental distribution of PrP may not correspond to a steady state of the system since the  
biological process is dramatically stopped with the animal death. Therefore the considerations about the stationary distribution might be not biologically relevant. Nevertheless we also observe bimodal distribution generically along the temporal dynamics,  a new feature in comparison to \cite{GreerPW}. 
For further modifications concerning model \eqref{eq:Greer} we suggest to take into account saturation in the fragmentation rate for large fibrils, as noticed in \cite{Caughey}.

Theoretical analysis of the model gives necessary criterion on the parameters for a bimodal distribution to occur in the stationary phase. Numerical analysis of  the eigenvalue problem shows in addition that parameters which lead to a bimodal distribution of PrPsc are not the most favourable for PrPsc accumulation in the exponential phase, see Figure \ref{fig:vp}. It suggests that infectious ability concentrated on a restricted range of PrPsc aggregate sizes (14 to 28 molecules \cite{Caughey}) does not consist in a replicative advantage for prion at the molecular level. Therefore, this restriction could be caused by a physico-chemical constraint, carried out by the host metabolism or by the prion aggregate itself. This could eventually be due to the PrPsc aggregate shape, since the most infectious intermediate aggregates have been found to be spherical or ellipsoidal whereas the others form fibrilar structures \cite{Caughey}. However, understanding the precise nature of the constraint would help to identify new therapeutic targets.

\medskip

Furthermore, we investigate the implications of varying the fragmentation and transconformation rates notably in the biological perspectives of the prion strain phenomenon.
Prion strains have been initially distinguished by incubation periods and lesion profiles in congenic mice \cite{Fraser, Bruce}.
Nowadays, a large body of literature suggests that differences between prion strains lie in the diversity of structures of PrPsc aggregates that can be stably and faithfully propagated  \cite{Cobb, Morales, Collinge, Safar, Legname, Aguzzi, Makarava, Thackray, Aguzzi08}. However, it remains poorly understood how these changes in the conformation of PrPsc aggregates can account for their physiopathological effects \cite{Legname}. Among the attempted biochemical characterisation of prion strains, a relationship was found between the relative stability values of PrPsc aggregates \cite{Legname} or level of aggregation \cite{Aguzzi} and incubation times,  indicating that less stable prions are more infectious, as judged by their shorter incubation times. This is presumably because unstable prions fragment more easily, giving rise to smaller aggregates of PrPsc that are more infectious than larger ones (represented by an increased transconformation rate in the mathematical model).
Our model agrees with these observations, since increasing the fragmentation rate $\beta $ leads to a faster PrPsc accumulation. More interestingly, our model enables to explore quickly the influence of every elementary mechanism  involved in prion replication.
Whereas changing parameters leads to similar effects on PrPsc accumulation kinetics, the resulting size-distribution variations are indeed different. The achievement of experimental size-distribution of PrPsc aggregates for many prion strains, as obtained by Silveira et al \cite{Caughey} for 263K prion strain, could therefore allow to better understand the mechanisms involved in prion strain phenomenon, by comparing the experimental distribution and the predicted ones. One limit to our approach is that prion strains are not only characterized by a precise incubation time but also by a specific cellular tropism. It would be of interest to introduce some mathematical formalism taking into account cell heterogenity in the brain and the resulting local influences.

In addition, transmission of prion diseases between different mammalian species is almost systematically less efficient than within a single species. This obstruction has been termed Species Barriers. Early studies argue that barrier resides in PrP primary structure difference between donor and recipient species  \cite{Prusiner90}. However, this issue has been called into question notably by BSE strain capability. Indeed, different strains propagated in the same host may thus have completely different barriers to another species. Consequently, transmission barrier appears to depend on prion strain specificities  \cite{Bruce94, Collinge95}. This is supported by our model. Indeed, we learn that the renormalized bimodal distribution of PrPsc aggregates tends to align along the eigenvector associated to the dominant eigenvalue $\Lambda(\overline V)$ (solid line in Figure \ref{fig:eigenvalue}(left)) \cite{Perthame_LN}. In our analysis, we assumed implicitely that the distribution of polymers is initially proportional to this eigenvector (which depends upon the the parameters of the model, and thus on the strain-host relationship). Otherwise there is some delay for the microscopic distribution to  reach the good shape which slows down the process of prion accumulation. This might account for  the strain adaptation mechanism, where primary inoculation is associated with more prolonged incubation period than subsequent passages in the same type of hosts \cite{Collinge}.    
According to these observations and the experimental ones, it should be very interesting to study experimentally the distribution of aggregates and infectivity between strains which are able or not to cross many species barriers. The distributions of PrPsc aggregates could be used to evaluate the influence of fragmentation, transconformation rate or infectivity distribution on the strenght of prion strain in crossing species barriers. These data could give a new signature of prion strain. Meanwhile the elucidation of structural basis of prion strains, strenght barrier species discrimination by these prion strain signature could be useful to predict the potential transmissibility of prion strain to human, notably for recent atypical scrapie strain, Bovine Amyloidic Spongiform Encephalopathy (BASE) strain and Chronic Wasting Disease strain \cite{Watts}.

Thus, our model could be used to investigate essential features of prion strains. The next step of our work is to approximate faithfully the inverse problem in order to obtain the size dependence of the transconformation rate from the  distribution of PrPsc. We therefore would be able to determine the most infectious compartment for several prion strains. This knowledge is a critical step for experimental approaches of prion infectivity investigation, like PMCA. PMCA goal is to quickly synthetize in vitro large amounts of PrPsc starting with minute
amounts of prions. The yield of this technique is very sensitive to experimental procedures. Then size distribution of infectivity could help to optimize strain specific PMCA protocols, notably by adaptating sonication steps to fall into the most infectious size compartment of PrPsc aggregates.

\medskip

Finally, we have analyse the stability of the steady states in a general framework. The difficulty arising here is that we cannot reduce the study to a set of Ordinary Differential Equations, as it is the case in \cite{GreerPW,GreDrieWanWeb}. Under general assumptions on the coefficients, compatible with bimodal distributions, the asymptotic stability of the healthy state (i.e. no prion aggregates) is established when the PrPc is low. We also prove that this healthy state is unstable when the PrPc production rate is high enough. This is in accordance with the asymptotic stability of the non-zero steady state for 'constant coefficients' proved in \cite{GreerPW}, even though a perfect dichotomy between the two results is left open in the general case. Biologically, these results can be interpreted as the propensity of PrPsc aggregates to give rise to prion disease depends on the amount of PrPc. In the more general context of understanding why some amyloids are infectious and others are not, it therefore could be useful to investigate the amount of amyloidogenic precursors.

\begin{acknowledgement*}
\noindent This work has been achieved during a visit of DO within the DEASE project (Marie Curie Early Stage Training multi Site (mEST) of the EU, MEST-CT-2005-021122). 
\\
We would like to sincerely thank  J. Silveira,  A. Hughson and B. Caughey for their experimental data about PrP distribution. We are also grateful to M. Doumic for fruitful discussions. This work would not have been achieved without the help of H. Zaag.
\end{acknowledgement*}


\end{document}